\documentclass[reqno,11pt]{amsart}
\usepackage{amssymb} 
\usepackage{latexsym} 
\usepackage{verbatim} 
\newif\ifpdf 
\ifpdf 
  \usepackage[pdftex]{graphicx} 
  \usepackage[pdftex]{hyperref} 
\else 
  \usepackage{graphicx} 
\fi 
\numberwithin{equation}{section} 
\theoremstyle{plain}

\theoremstyle{plain} 
\newtheorem{theo}{Theorem}[section] 
 
\newtheorem{cor}[theo]{Corollary} 
\newtheorem{defi}[theo]{Definition} 
\newtheorem{prop}[theo]{Proposition} 
\newtheorem{prop-def}[theo]{Proposition-definition} 
\newtheorem{lem}[theo]{Lemma} 
\newtheorem{remark}[theo]{Remark} 
\newtheorem{example}[theo]{Example} 
\newtheorem*{thmA}{Theorem A}  
\newtheorem*{thmB}{Theorem B}  
\newtheorem*{corC}{Corollary C}  
\newtheorem*{thmD}{Theorem D}  
\newtheorem*{corE}{Corollary E}  
\newtheorem*{thmF}{Theorem F}  
 
\newcommand{\CP}{\mathbf{P}} 
\newcommand{\er}{\mathbf{R}} 
\newcommand{\ku}{\mathbf{Q}} 
\newcommand{\co}{\mathbf{C}} 
\newcommand{\ze}{\mathbf{Z}}

\newcommand{\oh}{\mathcal{O}}

\newcommand{\fX}{\mathfrak{X}} 
 
\newcommand{\fD}{\mathfrak{D}} 
\newcommand{\ep}{\varepsilon} 
 
\newcommand{\al}{\alpha} 
 
\newcommand{\vol}{\operatorname{vol}} 
 
\newcommand{\NS}{N^1} 
\newcommand{\CNS}{CN^1} 
\newcommand{\Np}{N^p} 
\newcommand{\CNp}{CN^p} 
\newcommand{\cD}{\mathcal{D}} 
\begin{document} 
%
%
 
\setcounter{tocdepth}{2} 
 
\title[Differentiability of volumes of divisors]{Differentiability of volumes of divisors\\ and a problem of Teissier} 
\author{S{\'e}bastien Boucksom, Charles Favre, Mattias Jonsson} 
\address{CNRS-Universit{\'e} Paris 7\\ 
  Institut de Math{\'e}matiques\\ 
  F-75251 Paris Cedex 05\\ 
  France\\ 
  and\\ 
  Graduate School of Mathematical Sciences\\ 
  The University of Tokyo\\ 
  3-8-1 Komaba, Meguro-ku, Tokyo, 153-8914\\ 
  Japan}  
\email{boucksom@math.jussieu.fr} 
\address{Unidade Mista CNRS-IMPA\\ 
  Estrada Dona Castorina 110\\ 
  Rio de Janeiro / Brasil 22460-320\\ 
  and\\ 
  CNRS-Universit{\'e} Paris 7\\ 
  Institut de Math{\'e}matiques\\ 
  F-75251 Paris Cedex 05\\ 
  France} 
\email{favre@math.jussieu.fr} 
\address{Dept of Mathematics\\ 
  University of Michigan\\ 
  Ann Arbor, MI 48109-1043\\ 
  USA 
  \and 
  Dept of Mathematics\\ 
  KTH\\ 
  SE-100 44 Stockholm\\ 
  Sweden} 
\email{mattiasj@umich.edu,mattiasj@kth.se}  
\thanks{Second author 
  supported by the Japanese Society for the Promotion of Science. 
  Third author supported by NSF Grant No DMS-0449465, the  
  Swedish Science Council and the Gustafsson Foundation.} 
 
\date{\today} 
 
\begin{abstract} 
  We give an algebraic construction of the positive  
  intersection products of pseudo-effective classes first introduced 
  in~\cite{BDPP}, and use them to prove that the volume function  
  on the N{\'e}ron-Severi space of a projective variety is 
  $\mathcal{C}^1$-differentiable, expressing its differential as a 
  positive intersection product. We also relate the differential to 
  the restricted volumes 
  introduced in~\cite{ELMNP3,Takayama}. We then apply our 
  differentiability result to prove an algebro-geometric 
  version of the Diskant inequality in convex geometry, allowing us 
  to characterize the equality case of the  
  Khovanskii-Teissier inequalities for nef and big classes. 
\end{abstract}

 
\maketitle 
\setcounter{tocdepth}{1} 
\tableofcontents 
 
\section*{Introduction} 
The \emph{volume} of a line bundle $L$ on a projective variety 
$X$ of dimension $n$ is a nonnegative real number measuring the 
positivity of $L$ from the point of view of birational geometry. 
It is defined as the growth rate of sections of multiples of $L$: 
\begin{equation*} 
  \vol(L):=\limsup_{k\to\infty}\frac{n!}{k^n}\,h^0(X,kL) 
\end{equation*} 
and is positive iff the linear system $|kL|$ embeds $X$ 
birationally in a projective space for $k$ large enough; 
$L$ is then said to be \emph{big}. 
 
The volume  has been studied by several authors, and the 
general theory is presented in detail and with full references 
in~\cite[\S2.2.C]{Lazbible}. In particular, it is known that the volume only 
depends on the numerical class of $L$ in the real N{\'e}ron-Severi space 
$\NS(X)$, and that it uniquely extends to a continuous function on this whole 
space, such that $\vol^{1/n}$ is homogeneous of degree $1$,  concave 
on the open convex cone of big classes, and zero outside. 
 
Given its fundamental nature, it is quite natural to ask what kind of 
regularity besides continuity the volume function exhibits in general. 
In the nice survey~\cite{ELMNP1}, many specific examples were 
investigated, leading the authors to conjecture that the volume 
function is always real analytic on a "large" open subset of the big 
cone. Our concern here will be the differentiability of the volume 
function. The simple example of $\CP^2$ blown-up in one point already 
shows that the volume function is not twice differentiable on the 
entire big cone in general. Our main result is 
\begin{thmA} 
  The volume function is $\mathcal{C}^1$-differentiable on the big 
  cone of $\NS(X)$. If $\alpha\in \NS(X)$ is big and  
  $\gamma\in \NS(X)$ is arbitrary, then 
  \begin{equation*} 
    \left. \frac{d}{dt}\right|_{t=0} 
    \vol(\alpha+t\gamma)=n\langle\alpha^{n-1}\rangle\cdot\gamma. 
  \end{equation*} 
\end{thmA} 
The right-hand side of the equation above involves the \emph{positive 
intersection product} $\langle\alpha^{n-1}\rangle\in \NS(X)^{\ast}$ of 
the big class $\alpha$, first introduced in the analytic context 
in~\cite{BDPP}. We shall return to its algebraic definition later in 
this introduction, when discussing our method of proof. 
 
\smallskip 
We then proceed to show that the derivative of the volume in the 
direction of a class determined by a prime divisor can also be 
interpreted as a restricted volume, as introduced and studied in 
\cite{ELMNP3}. Recall that if $V$ is a subvariety of $X$, the 
\emph{restricted volume} on $V$ of a line bundle $L$ on $X$ measures 
the growth of sections in $H^0(V,kL|_{V})$ that extend to $X$. It is 
defined as 
\begin{equation*} 
  \vol_{X|V}(L):=\limsup_{k\to\infty}\frac{d!}{k^d}~h^0(X|V,kL) 
\end{equation*} 
where $d:=\dim V$ and $h^0(X|V,kL)$ denotes the rank of the restriction map 
\begin{equation*} 
  H^0(X,kL)\to H^0(V,kL|_{V}). 
\end{equation*} 
Restricted volumes have recently 
played a crucial role in the proof of the boundedness of 
pluricanonical systems of varieties of general type in 
\cite{Takayama} and implicitly in~\cite{hacon}. Here we relate the 
positive intersection products with restricted volumes, and show 
\begin{thmB} 
  If $D$ is a prime divisor on the smooth projective variety $X$  
  and $L$ is a big line bundle 
  on $X$, then the restricted volume of $L$ on $D$ satisfies 
  \begin{equation*} 
    \vol_{X|D}(L)=\langle L^{n-1}\rangle\cdot D. 
  \end{equation*} 
\end{thmB} 
This statement in particular implies that the restricted volume only 
depends on the numerical class of both $L$ and $D$ in $\NS(X)$.  When 
$V$ is an irreducible subvariety,~\cite{ELMNP3} and~\cite{Takayama} 
have independently shown that the restricted volume $\vol_{X|V}(L)$ of 
a big line bundle $L$ can be expressed as the asymptotic intersection 
number of the moving parts of $|kL|$ with the strict transforms of $V$ 
on appropriate birational models $X_k$ of $X$, when $L$ satisfies an 
additional positivity assumption along $V$.  The main (and difficult) 
result of~\cite{ELMNP3} says that $\vol_{X|V}(L)=0$ otherwise. Our 
contribution in Theorem~B is to show that the asymptotic intersection 
number along a \emph{prime divisor} $V=D$ in fact coincides with the 
intersection number $\langle L^{n-1}\rangle\cdot D$ when the 
additional positivity condition is satisfied, using our 
differentiability result, and that both sides are $0$ otherwise, 
relying on the orthogonality of Zariski decompositions instead of the 
main result of~\cite{ELMNP3}. 
 
Theorems~A and~B yield the following corollary, 
which was kindly communicated to us by 
R.~Lazarsfeld and M.~Musta\c{t}\v{a} 
and which 
inspired our differentiability result. 
\begin{corC} 
  If $D$ is a prime divisor on the smooth projective variety $X$ and  
  $L$ is a big line bundle, then 
  \begin{equation*} 
    \left. \frac{d}{dt}\right|_{t=0}\vol(L+tD)=n\vol_{X|D}(L). 
  \end{equation*} 
\end{corC} 
We also give an application of our differentiability theorem and 
characterize the equality case in the Khovanskii-Teissier inequalities 
for big and nef classes, a problem considered by Teissier in 
\cite[p.96]{Teissier82} and~\cite[p.139]{Teissier88}.  Recall that one 
version of the Khovanskii-Teissier inequalities~\cite{Teissier79} for 
a pair of nef classes $\alpha,\beta\in \NS(X)$ asserts that the 
sequence $ k \mapsto \log ( \alpha^k \cdot \beta^{n-k})$ is 
concave. Here we prove 
\begin{thmD} 
  If $\alpha,\beta\in \NS(X)$ are big and nef classes, 
  then the following are equivalent: 
  \begin{itemize} 
  \item[(i)] 
    the concave sequence $k\mapsto\log(\alpha^k\cdot\beta^{n-k})$ 
    is affine; 
  \item[(ii)] 
    $(\alpha^{n-1}\cdot\beta)=(\alpha^n)^{1-\frac1n}(\beta^n)^{\frac1n}$; 
  \item[(iii)] 
    $\alpha$ and $\beta$ are proportional. 
  \end{itemize} 
\end{thmD} 
An equivalent formulation is 
\begin{corE} 
  The function $\alpha\mapsto(\alpha^n)^{\frac 1 n}$ is strictly 
  concave on the big and nef cone of $\NS(X)$. 
\end{corE} 
Note that when $\alpha$ and $\beta$ are ample, these results are 
direct consequences of the Hodge-Riemann relations, see for 
instance~\cite[Theorem~6.32]{voisin}.  However, our approach is 
different and purely algebraic.  In fact, the two statements above and 
their proofs are inspired by their counterparts in convex 
geometry~\cite{Schneider}.  When $X$ is a projective toric variety, a 
big and nef class $\alpha\in \NS(X)$ corresponds to a convex polytope 
in $\er^n$ (unique up to translation) with Euclidean volume 
$(\alpha^n)/{n!}>0$.  The Khovanskii-Teissier inequalities are known 
in this setting as the Alexandrov-Fenchel inequalities for mixed 
volumes, a far-reaching generalization of the classical isoperimetric 
inequality.  A different, ``effective'', strengthening of the 
isoperimetric inequality was given by Bonnesen~\cite{bonnesen} (in 
dimension two) and Diskant~\cite{Diskant} (in any dimension).  It 
immediately implies the analog of Theorem~D in convex geometry; in 
particular, equality holds in the isoperimetric inequality iff the 
convex body is a ball.  The proof by Diskant is based on the 
differentiability of the volume function of (inner parallel) convex 
bodies, a fact established by Alexandrov.  Following the same 
strategy, and using Theorem~A, we prove the following version of the 
Diskant inequality in our setting, thus providing a solution 
to~\cite[Problem~B]{Teissier82}. 
\begin{thmF} 
  If $\alpha,\beta\in \NS(X)$ are big and nef classes and $s$ is the 
  largest real number such that $\alpha- s\beta$ is pseudo-effective, 
  then 
  \begin{equation*} 
    (\alpha^{n-1}\cdot\beta)^{\frac{n}{n-1}}-(\alpha^n) 
    (\beta^n)^{\frac{1}{n-1}} \geq  
    \left( 
      (\alpha^{n-1}\cdot\beta)^{\frac{1}{n-1}}-s(\beta^n)^{\frac{1}{n-1}} 
    \right)^n. 
  \end{equation*} 
\end{thmF} 
Theorem~D immediately follows from Theorem~F. 
 
\medskip 
Let us now present our method of proof of Theorem~A. 
It is a fundamental fact that the volume of an 
arbitrary big line bundle $L$ admits an 
intersection-theoretic interpretation, generalizing the equality 
$\vol(L)=(L^n)$ when $L$ is ample. Indeed, a remarkable theorem of 
Fujita~\cite{Fujita} (see also~\cite{DEL}) states that $\vol(L)$ is 
the supremum of all intersection numbers $(A^n)$, where $A$ is an 
ample $\ku$-divisor and $E$ is an effective $\ku$-divisor on a 
modification $\pi:X_\pi\to X$ such that $\pi^{\ast}L=A+E$. 
 
In view of this result, it is natural to put \emph{all} birational 
models of $X$ on equal footing, and study numerical classes defined  
on all birational models at the same time.  
Being purely algebraic, this technique 
extends to any projective variety 
over any algebraically closed field of characteristic zero. 
 
More precisely, we introduce the \emph{Riemann-Zariski space} $\fX$ of 
$X$ as the projective limit of all birational models of $X$. A 
(codimension $p$ numerical) class in $\fX$ is then a collection of 
(codimension $p$ numerical) classes in each birational model of $X$ 
that are compatible under push-forward.  The set of all these classes 
is an infinite dimensional vector space that we denote by $\Np(\fX)$. 
It naturally contains the space $\CNp(\fX)$, the union of the spaces 
$\Np(X_\pi)$ of numerical classes of codimension $p$ cycles of all 
birational models $X_\pi$ of $X$. One can extend to $\CNS(\fX)$ the 
usual notions of pseudo-effective, big and nef classes.  
Note that related 
objects have already been introduced in the context of Mori's minimal 
model program by Shokurov~\cite{Shokurov}. The notion of 
\emph{b-divisors}, crucial to his approach, can be interpreted as 
divisors on $\fX$. Chow groups on the 
Riemann-Zariski space also appeared in the work of 
Aluffi~\cite{aluffi} and numerical  
classes in~\cite{BFJ1,cantat,manin} in the case of surfaces. 
See also~\cite{valtree,BFJ2} for a local study. 
 
We then introduce the notion of  
\emph{positive intersection product} 
\begin{equation*} 
  \langle\alpha_1\cdot\ldots\cdot\alpha_p\rangle\in \Np(\fX) 
\end{equation*} 
of big 
classes $\alpha_1,\dots,\alpha_p\in \CNS(\fX)$, defined as the least 
upper bound of the products $\beta_1\cdot\ldots\cdot\beta_p$ of nef 
classes $\beta_i$ on modifications $\pi:X_\pi\to X$ such that 
$\beta_i\leq\pi^{\ast}\alpha_i$. Although this product is non-linear, 
it is homogeneous and increasing in each variable, and continuous on the 
big cone.  In this terminology, Fujita's theorem is equivalent to the 
relation 
\begin{equation*} 
  \vol(\alpha)=\langle\alpha^n\rangle 
\end{equation*} 
for any big class $\alpha\in\NS(X)$.  
Using the easy but fundamental inequality 
\begin{equation*} 
  \vol(A-B)\geq(A^n)-n(A^{n-1}\cdot B) 
\end{equation*} 
for any two nef Cartier classes $A,B\in\CNS(\fX)$, we deduce a 
sub-linear control for the volume function  
$\vol(\alpha+t\gamma) 
\geq\vol(\alpha)+tn\langle\alpha^{n-1}\rangle\cdot\gamma+O(t^2)$ for 
any two Cartier classes $\alpha,\gamma\in\CNS(\fX)$ such that $\alpha$ 
is big, from which we easily infer Theorem A.  As an immediate 
consequence, we get the following orthogonality relation 
\begin{equation*} 
  \langle\alpha^n\rangle=\langle\alpha^{n-1}\rangle\cdot\alpha 
\end{equation*} 
for any psef class $\alpha\in\CNS(\fX)$, which was the crucial point 
in the dual characterization of pseudo-effectivity of~\cite{BDPP}. 
 
\medskip 
In a similar vein, the proof of Theorem~B is based on 
a suitable generalization of Fujita's theorem for restricted volumes, 
obtained independently in~\cite{ELMNP3} and~\cite{Takayama}.  If 
$D$ is a prime divisor on a smooth projective variety $X$, 
we say that a line bundle $L$ is \emph{$D$-big} if there exists a 
decomposition $L=A+E$, where $A$ is an ample 
$\ku$-divisor and $E$ is an effective $\ku$-divisor on $X$ whose 
support does not contain $D$. The generalization of Fujita's theorem 
expresses the restricted volume of a $D$-big line bundle as the 
supremum of all intersection numbers $A^{n-1}\cdot D_\pi$, where 
$\pi:X_\pi\to X$ is a modification, $D_\pi$ denotes the strict 
transform of $D$, $A$ is an ample $\ku$-divisor and $E$ is an 
effective $\ku$-divisor on $X_\pi$ whose support does not contain 
$D_\pi$ and such that $\pi^{\ast}L=A+E$.  In Section~\ref{sec:res}, 
we again interpret this result in the framework of the Riemann-Zariski 
space and define the restricted positive intersection product 
\begin{equation*} 
  \langle\alpha_1\cdot\ldots\cdot\alpha_p\rangle|_{\fD} 
\end{equation*} 
of $D$-big classes $\alpha_i$ on $X$ as a numerical class of the 
Riemann-Zariski space $\fD$ associated to $D$, in such a way that the 
restricted positive intersection product $\langle 
L^{n-1}\rangle|_{\fD}$ coincides with the asymptotic intersection 
number $||L^{n-1}\cdot D||$ of $L$ and $D$ as introduced 
in~\cite{ELMNP3} if $L$ is $D$-big. The generalized Fujita theorem 
mentioned above then reads 
\begin{equation*} 
  \vol_{X|D}(L)=\langle L^{n-1}\rangle|_{\fD} 
\end{equation*} 
when $L$ is $D$-big, and Theorem~B asserts that the restricted 
positive intersection product $\langle L^{n-1}\rangle |_{\fD}$ 
coincides with the intersection number  
$\langle L^{n-1}\rangle\cdot D$. 
 
The former is computed as a limit as $k\to\infty$ of intersection 
numbers of the form $A_k^{n-1}\cdot D_k$ with $A_k$ an ample 
$\ku$-divisor $\le\pi_k^*L$ on a blow-up $\pi_k:X_k\to X$ and $D_k$ 
the strict transform of $D$ on $X_k$, and the latter is then the limit 
of $A_k^{n-1}\cdot\pi_k^*D$.  We easily deduce from this relation that 
$\langle L^{n-1}\rangle |_{\fD}\leq\langle L^{n-1}\rangle\cdot D$. 
The reverse inequality is obtained by showing that 
\begin{equation*} 
  \limsup_{t\to 0+}\frac1t\left(\vol(L)-\vol(L-tD)\right) 
  \leq  n\vol_{X|D}(L), 
\end{equation*} 
building on the basic relation $h^0(X,L)=h^0(X,L-D)+h^0(X|D,L)$.  This 
yields Theorem~B in the $D$-big case.  When $L$ is big but not 
$D$-big, we prove that the inequality $\langle L^{n-1}\rangle\cdot 
D\geq\vol_{X|D}(L)$ still holds, and that $\langle L^{n-1}\rangle 
\cdot D=0$, as a consequence of the orthogonality relation mentioned 
above. We get in particular that $\vol_{X|D}(L)=0$ when $L$ is not 
$D$-big, thus reproving by a different method a special case of the 
main (and difficult) result of~\cite{ELMNP3}. 
 
\medskip 
Let us briefly indicate the structure of our article. We introduce 
numerical classes in the Riemann-Zariski space in 
Section~\ref{sec:rie}, and study their positivity properties. The 
positive intersection product is then defined in 
Section~\ref{sec:pos}. We turn to the volume function on the set 
of psef classes in the Riemann-Zariski space in Section~\ref{sec:vol}. 
In particular, this section contains the proof of the 
Differentiability Theorem, and its applications to the 
characterization of the equality case in the Khovanskii-Teissier 
inequalities: Theorems~D and~F and 
Corollary~E. We also include in Section~\ref{sec:zar} an 
informal discussion making a link between  
positive intersection products and 
Zariski-type decompositions for psef classes, intended to shed further 
light on our construction.  The paper concludes with  
a presentation of the restricted volume from the  
Riemann-Zariski point of view in Section~\ref{sec:res},  
allowing us to prove Theorem~B and Corollary~C. 
 
\medskip 
\noindent\textbf{Acknowledgment}.  It is a great pleasure to thank the 
authors of~\cite{ELMNP3} for sharing their work with us.  Their ideas 
have obviously had a great influence on the present paper. In 
particular, the statement of Corollary~C was suggested to us by 
R.~Lazarsfeld and M.~Musta\c{t}\v{a}, and they informed us after 
completion of the present paper that they had independently obtained a 
proof of it.  We also thank B.~Teissier for explaining the history of 
the problems considered here; S.~Takayama for interesting discussions 
and for sharing his preprint; and F.~Wittke for spotting a mistake in 
a previous version of the paper. 
 
 
\section{Classes on the Riemann-Zariski space}\label{sec:rie} 
 
 
\subsection{The Riemann-Zariski space} 
Let $X$ be a projective variety. Since we shall deal with 
classes living on arbitrary modifications of $X$, it is convenient to 
introduce the following terminology. By a \emph{blow-up} of $X$, we 
mean a birational morphism $\pi:X_\pi\to X$, where $X_\pi$ is a 
\emph{smooth} projective variety. If $\pi$ and $\pi'$ are two blow-ups 
of $X$, we say that $\pi'$ \emph{dominates} $\pi$, and write 
$\pi'\geq\pi$, if there exists a birational morphism (necessarily 
unique) $\mu:X_{\pi'}\to X_\pi$ such that $\pi'=\pi\circ\mu$. This 
endows the set of blow-ups of $X$ with a 
partial order relation. By the Hironaka desingularization theorem, 
this ordered set is nonempty and forms a directed family.  The 
\emph{Riemann-Zariski space} of $X$ is the projective limit 
\begin{equation*} 
  \fX:=\varprojlim_\pi X_\pi. 
\end{equation*} 
We refer to~\cite{ZS,vaquie} for a thorough 
discussion of the structure of this space, which is 
introduced here merely in order to make later definitions more suggestive.

 
\subsection{Weil and Cartier classes} 
For any smooth projective variety $Y$ of dimension $n$, and any 
integer $0\le p\le n$, let $\Np(Y)$ be the real vector space of 
numerical equivalence classes of codimension $p$ cycles.  It is a 
finite dimensional space.  Any birational morphism $\mu:Y'\to Y$ 
induces in a contravariant way a pull-back morphism 
\begin{equation*} 
  \mu^{\ast}:\Np(Y)\to \Np(Y')~; 
\end{equation*} 
and in a covariant way a push-forward morphism 
\begin{equation*} 
  \mu_{\ast}:\Np(Y')\to \Np(Y). 
\end{equation*} 
There is an intersection pairing 
$\Np(Y)\times N^{n-p}(Y)\to \er$, which is 
preserved under pull-back by birational morphisms,  
and for which push-forward and pull-back are 
adjoint to each other. We refer to~\cite[Chapter~19]{Fulton} for further 
details on the space $N^p$ and its link with the cohomology groups 
$H^{p,p}$ in the case of a complex projective variety. 
 
We now consider an arbitrary projective variety $X$.  The inverse 
family of blow-ups $\pi:X_\pi\to X$ induces two families of arrows 
between the spaces $\Np(X_\pi)$: one is the inverse family with 
push-forward morphisms as arrows  
$\mu_{\ast}: \Np(X_{\pi'}) \to\Np(X_{\pi})$,  
whenever $\mu:X_{\pi'}\to X_\pi$, and the other is the 
direct family with pull-back morphisms as arrows  
$\mu^{\ast}:\Np(X_{\pi}) \to \Np(X_{\pi'})$. 
 
\medskip 
\begin{defi} {~} 
  \begin{itemize} 
  \item The space of $p$-codimensional Weil classes on the 
    Riemann-Zariski space $\fX$ is defined as the projective limit 
    \begin{equation*} 
      \Np(\fX):=\varprojlim_\pi \Np(X_\pi) 
    \end{equation*} 
    with arrows defined 
    by push-forward. It is endowed with its projective limit topology, 
    which will be called the \emph{weak topology}. 
  \item The space of $p$-codimensional Cartier classes on $\fX$ is 
    defined as the inductive limit 
    \begin{equation*} 
      \CNp(\fX):=\varinjlim_\pi \Np(X_\pi) 
    \end{equation*} 
    with arrows induced by 
    pull-back. It is endowed with its inductive limit topology, 
    which will be called the \emph{strong topology}. 
  \end{itemize} 
\end{defi} 
The terminology stems from the fact that Weil divisors are meant to be 
pushed-forward, whereas Cartier divisors are meant to be pulled-back. 
Note that these spaces are infinite dimensional as soon as  
$\dim X\geq 2$ and $p \neq 0,n$. 
 
A Weil class $\alpha\in \Np(\fX)$ is by definition described by its 
\emph{incarnations} $\alpha_\pi\in\Np(X_\pi)$ on each smooth 
birational model $X_\pi$ of $X$, compatible with each other by 
push-forward.  
Convergence of a sequence (or a net) $\al_n \to \al$ in the  
weak topology means 
$\alpha_{n,\pi} \to \alpha_\pi$ in $\Np(X_\pi)$ for each 
$\pi$. 
 
On the other hand, the relation $\mu_{\ast}\mu^{\ast}\alpha=\alpha$ 
when $\mu$ is a birational morphism shows that there is an injection 
\begin{equation*} 
  \CNp(\fX)\to \Np(\fX) 
\end{equation*} 
i.e.\ a Cartier class is in particular a Weil 
class.  Concretely, a Weil class $\alpha$ is Cartier iff there exists 
$\pi$ such that its incarnations $\alpha_{\pi'}$ on higher blow-ups 
$X_{\pi'}$ are obtained by pulling-back $\alpha_\pi$. We will call 
such a $\pi$ a \emph{determination} of $\alpha$. Conversely, there are 
natural injective maps 
\begin{equation*} 
  \Np(X_\pi)\hookrightarrow \CNp(\fX) 
\end{equation*} 
which extend a given class 
$\beta\in\Np(X_\pi)$ to a Cartier class by pulling it back, and of 
course this Cartier class is by construction determined by $\pi$. In 
the sequel, we shall always identify a class $\beta \in \Np(X_\pi)$ 
with its image in $\CNp(\fX)$. 
 
Note that the topology induced on $\CNp(\fX)$ by the weak topology of 
$\Np(\fX)$ is coarser than the strong topology. By definition, a map 
on $\CNp(\fX)$ is continuous in the 
strong topology iff its restriction to each of the  
finite dimensional subspaces $\Np(X_\pi)$ is continuous.   
 
Finally we note that 
\begin{lem}\label{dense} 
  The natural continuous injective map 
  \begin{equation*} 
    \CNp(\fX)\to \Np(\fX) 
  \end{equation*} 
  has dense image (in the weak topology). 
\end{lem} 
\begin{proof} 
  If $\alpha\in \Np(\fX)$ is a given Weil class, then we 
  can consider the Cartier classes $\alpha_\pi$  
  determined by the incarnation of $\alpha$ on 
  $X_\pi$, and it is obvious that the net $\alpha_\pi$ converges to 
  $\alpha$ in the weak topology as $\pi\to\infty$, since 
  $\alpha_\pi$ and $\alpha$ have by definition the same 
  incarnation on $X_\pi$. 
\end{proof} 
 
 
\subsection{Divisors} 
We shall refer to the space $\CNS(\fX):=\varinjlim_\pi \NS(X_\pi)$ as the 
\emph{N{\'e}ron-Severi space} of $\fX$. Its elements are related to the 
so-called b-Cartier divisors on $X$, in the sense of 
Shokurov~\cite{Shokurov,b-div}.  These are by definition elements of 
the inductive limit 
\begin{equation*} 
  \mathrm{CDiv}(\fX):=\varinjlim_\pi\mathrm{Div}(X_\pi) 
\end{equation*} 
of the 
spaces of (Cartier) $\er$-divisors on each $X_\pi$. According to our 
present point of view, we will change the terminology and call an 
element of $\mathrm{CDiv}(\fX)$ a \emph{Cartier divisor on} $\fX$. 
 
It is then immediate to see that the N{\'e}ron-Severi space $\CNS(\fX)$ is 
the quotient of $\mathrm{CDiv}(\fX)$ modulo numerical equivalence. In 
particular, a Cartier divisor $D$ on $X$ itself, however singular this 
space, induces a Cartier class in $\CNS(\fX)$. 
 
 
\subsection{Positive classes} 
When $Y$ is a smooth projective variety, we define the  
\emph{psef cone} (a short-hand for pseudo-effective cone)  
of $\Np(Y)$ as the 
closure of the convex cone generated by effective codimension $p$ 
cycles. 
\begin{prop} 
  For any smooth projective variety $Y$,  
  $\Np(Y)$ is a finite dimensional 
  real vector space in which the psef cone is  
  convex, closed, strict (i.e.\ $\pm\alpha$ psef implies $\alpha=0$) 
  and has a compact basis (i.e.\ the set of classes 
  $\beta\in \Np(Y)$ with $\alpha -\beta$ psef is compact for every 
  psef class $\alpha$). 
\end{prop} 
\begin{proof} 
  The fact that $\Np(Y)$ is finite dimensional is proved 
  in~\cite[Example 19.1.4]{Fulton}. By definition, the psef cone is 
  convex and closed. Any strict closed convex cone in a finite 
  dimensional space has a compact basis, so it remains to prove that 
  the psef cone is strict. If $\alpha\in N^p(Y)$ is psef, then clearly 
  $\alpha\cdot h_1\cdot\ldots\cdot h_{n-p}\geq 0$ for all ample 
  classes $h_i$ on $Y$. Therefore if $\pm\alpha$ are psef, we get that 
  $\alpha$ is zero against any complete intersection class, 
  i.e. $\alpha\cdot\beta_1\cdot\ldots\cdot\beta_{n-p}=0$ for all 
  classes $\beta_i\in N^1(Y)$ (since any such class can be written as 
  a difference of ample classes).  Now if $f:Z\to Y$ is a surjective 
  smooth morphism, $\pm f^{\ast}\alpha$ is also psef, and thus also 
  zero against any complete intersection class on $Z$. Applying this 
  to $f:Z=\mathbf{P}(E)\to Y$ for a vector bundle $E$ on $Y$ yields 
  that $\alpha$ is zero against the $(n-p)$-th Segre class of any 
  vector bundle $E$ on $Y$, which is the push-forward of the 
  appropriate power of the tautological line bundle of 
  $\mathbf{P}(E)$. Since such Segre classes generate $N^{n-p}(Y)$, 
  this finally shows that $\alpha$ is zero in $N^p(Y)$. 
 \end{proof} 
In the sequel, we shall write $\alpha\geq 0$ when $\alpha\in\Np(Y)$  
is a psef class. 
If $\mu:Y'\to Y$ is a birational morphism, the push forward 
$\mu_{\ast}\alpha\in\Np(Y)$ of a psef class $\alpha\in 
\Np(Y')$ is also psef, thus we can introduce the following 
definition. 
\begin{defi} 
  A Weil class $\al \in \Np(\fX)$ is \emph{psef} iff all 
  its incarnations $\al_\pi \in \Np(X_\pi)$ are psef. 
\end{defi} 
The set of all psef classes is a strict convex cone in $\Np(\fX)$ that 
is closed in the weak topology.  We will write $\alpha\geq 0$ if 
$\alpha\in \Np(\fX)$ is psef. 
\begin{prop}\label{compact}  
  The psef cone has compact basis. In other   words, 
  if $\alpha\in \Np(\fX)$ is a given psef class, the set of Weil 
  classes $\beta\in \Np(\fX)$ such that $0\leq\beta\leq\alpha$ is 
  compact (in the weak topology). 
\end{prop} 
\begin{proof} 
  The set in question is the projective limit of the 
  corresponding sets  
  $K_\pi:=\{\gamma\in\Np(X_\pi)\ |\ 0\leq\beta_\pi\leq\alpha_\pi\}$, 
  and each of these sets is compact 
  as recalled above, so the result follows from the Tychonoff theorem. 
\end{proof} 
We now consider positive Cartier classes on $\fX$.  
If $Y$ is a projective variety, recall~\cite[\S1.4.C]{Lazbible}  
that the \emph{nef cone} of $\NS(Y)$ is the 
closure of the cone of classes determined by ample divisors.  A nef 
class is thus psef.  The interior of the psef cone of  
$N^1(Y)$ is called the 
\emph{big cone}.  If $\mu:Y'\to Y$ is a birational morphism, then a 
class $\alpha\in \NS(Y)$ is nef (resp.\ psef, big) iff 
$\mu^{\ast}\alpha$ is. This property enables us to extend 
the definitions to the Riemann-Zariski space. 
\begin{defi} 
  A Cartier class $\alpha\in \CNS(\fX)$ is said to 
  be nef (resp.\ psef, big) iff its incarnation $\alpha_\pi$ is nef 
  (resp.\ psef, big) for some determination $\pi$ of $\alpha$. 
\end{defi} 
The set of psef classes of $\CNS(\fX)$ will be called the \emph{psef 
  cone} of $\CNS(\fX)$. It is obviously a closed convex cone, and in 
fact it coincides with the inverse image of the psef cone of 
$\NS(\fX)$ under the continuous injection $\CNS(\fX)\to \NS(\fX)$.  In 
other words, a Cartier class $\alpha\in \CNS(\fX)$ is psef as a 
Cartier class iff it is psef as a Weil class, thus the terminology 
chosen is not ambiguous. 
\begin{prop}\label{prop:recall} 
  The psef cone of $\CNS(\fX)$ is the 
  closure of the big cone, i.e.\ the set of big classes. The big cone 
  is the interior of the psef cone. 
\end{prop} 
\begin{proof}  
  The analogous statement is true in each space 
  $\NS(X_\pi)$, so the result follows easily from the 
  behavior of psef and big classes under pull-back recalled above and 
  the definition of the strong topology. 
\end{proof} 
The set of nef classes of $\CNS(\fX)$ will be called the  
\emph{nef cone} of $\CNS(\fX)$. It is a closed convex cone. 
\begin{remark} 
  A b-divisor is b-big (resp.\ b-nef) 
  in the terminology of Shoku\-rov iff its class in $\CNS(\fX)$ 
  is big (resp.\ nef). 
\end{remark} 
\begin{remark} 
  There is no "ample cone" on $\fX$, i.e.\ the interior of the nef 
  cone of $\CNS(\fX)$ is empty. Indeed, if $\alpha$ is a nef class 
  determined by $\pi$, then for any strictly higher blow-up 
  $\pi'\geq\pi$ its incarnation $\alpha_{\pi'}\in 
  \NS(X_{\pi'})$ lies on the boundary of the nef cone. 
\end{remark} 
\begin{remark}\label{rem:nef} 
  One can also define the nef cone in $\NS(\fX)$ as the 
  closure in the weak topology of the  nef cone of  
  $\CNS(\fX)$: see~\cite{BFJ2} for details in a similar 
  situation.  
\end{remark} 
 
 
\subsection{Toric varieties}\label{sec:toric} 
Let $X$ be an $n$-dimensional projective toric variety, 
see~\cite{odabook} for background. If we restrict to toric blow-ups 
$\pi:X_\pi\to X$ in the definitions above, we obtain a toric 
Riemann-Zariski space $\fX_\mathrm{tor}$ and a toric N{\'e}ron-Severi 
space $\CNS(\fX_{\mathrm{tor}})\subset \CNS(\fX)$. One easily checks 
that $\CNS(\fX_{\mathrm{tor}})$  
is canonically isomorphic to the space of functions 
$g:N_\er\to\er$ that are piecewise linear with respect to some 
rational fan decomposition of $N_\er$, modulo linear forms.  Here as 
usual $N_\er=N\otimes_\mathbf{Z}\er$, where $N=\ze^n$ is the lattice of 1-parameters subgroups of the torus $(\co^{\ast})^n$. 
 
To a piecewise linear function $g$ as above is associated its Newton 
polytope $\mathrm{Nw}(g)$ in the dual space $M_\er = N_\er^*$. It is 
defined as the set of linear forms $m\in M_\er$ such that $m\leq g$ 
everywhere on $N_\er$. This polytope only depends on the cohomology 
class $\alpha\in \CNS(\fX_{\mathrm{tor}})$, up to translation, thus we 
will denote it by $\mathrm{Nw}(\alpha)=\mathrm{Nw}(g)$. It is a 
standard fact in toric geometry that: 
\begin{itemize} 
\item[(i)] $\alpha$ is psef iff $\mathrm{Nw}(\alpha)$ is nonempty; 
\item[(ii)] $\alpha$ is big iff $\mathrm{Nw}(\alpha)$ has nonempty interior; 
\item[(iii)] $\alpha$ is nef iff $g$ is convex. 
\end{itemize} 
Note that when $\alpha$ is nef, $g$ is equal to the maximum of 
all linear forms lying in $Nw(\alpha)$.  A toric nef class can hence 
be recovered from its Newton polytope, and we can identify toric nef 
classes and polytopes in $M_\er=\er^n$ with rational maximal faces, up 
to translation. 
 
 
\section{Positive intersection product}\label{sec:pos} 
 
\subsection{Intersections of Cartier classes} 
On a smooth projective variety $Y$, the intersection product of $p$ 
classes $\alpha_1,\dots,\alpha_p\in \NS(Y)$ defines an element 
$\alpha_1\cdot\ldots\cdot\alpha_p\in \Np(Y)$. Recall that for any 
birational morphism $\mu: Y' \to Y$, we have 
$\mu^*\alpha_1\cdot\ldots\cdot\mu^*\alpha_p= 
\mu^*(\alpha_1\cdot\ldots\cdot\alpha_p)$, see~\cite[Chapter~19]{Fulton}. 
 
We now define the intersection product of $p$ Cartier classes 
$\alpha_1,\dots,\alpha_p\in \CNS(\fX)$ having a common determination on 
$X_\pi$ as the Cartier class in $\CNp(\fX)$ determined by 
$\alpha_{1,\pi}\cdot\ldots\cdot\alpha_{p,\pi}$. Since the intersection 
product is preserved by pull-back, this is independent on the choice 
of determination. It is moreover continuous for the strong topology on 
$\CNS(\fX)$. 
\begin{remark}\label{rem:extend} 
  By the projection formula we have  
  $f_*( f^*\alpha\cdot\beta)=\alpha\cdot f_*\beta$  
  for any birational morphism $f$ between 
  smooth projective varieties.  
  It is thus possible to extend the 
  intersection product to a bilinear map 
  $\CNp(\fX)\times\mathrm{N}^{q}(\fX)\to\mathrm{N}^{p+q}(\fX)$. 
  In the case $p=1$ and $q=n-1$, we get a 
  pairing $\CNS(\fX) \times \mathrm{N}^{n-1}(\fX)\to\mathbf{R}$, and 
  under this pairing the psef cone of $\mathrm{N}^{n-1}(\fX)$ is 
  the dual of the nef cone in $\CNS(\fX)$. 
\end{remark} 
 
When the $\alpha_i$ are nef classes, it is easy to see that 
$\alpha_1\cdot\ldots\cdot\alpha_p\in \CNp(\fX) \subset \Np(\fX)$ is psef. 
More generally, we have 
\begin{lem} 
  If $\alpha_i\in \CNS(\fX)$, $1\le i\le p$ are Cartier classes 
  with $\alpha_1$ psef and $\alpha_i$ nef for $i\geq 2$, then their 
  intersection product $\alpha_1\cdot\ldots\cdot\alpha_p\in \Np(\fX)$ 
  is psef. 
\end{lem} 
\begin{proof} 
  By continuity, it is enough to check this when for some 
  common determination $\pi$, the class $\alpha_1$ is represented on 
  $X_\pi$ by a effective divisor and the $\alpha_i$ are represented by 
  very ample divisors for $i\geq 2$. But then the result is clear. 
\end{proof} 
The lemma implies the following monotonicity property which 
will be crucial in what follows. 
\begin{prop}\label{p:increase} 
  Let $\alpha_i$ and $\alpha_i'$, $1\le i\le p$, 
  be nef classes in $\CNS(\fX)$,  
  and suppose that  $\alpha_i\geq\alpha_i'$ for   $i=1,\dots,p$.  
  Then we have 
  \begin{equation*} 
    \alpha_1\cdot\ldots\cdot\alpha_p\geq\alpha_1'\cdot\ldots\cdot\alpha_p' 
  \end{equation*} 
  in $\Np(\fX)$. 
\end{prop} 
\begin{proof} 
  By successively replacing each $\alpha_i$ by $\alpha_i'$ 
  and using the symmetry of the intersection product, it is enough to 
  consider the case where $\alpha_i=\alpha_i'$ for $i\geq 2$. But then 
  $(\alpha_1-\alpha_1')\cdot\alpha_2\cdot \ldots\cdot\alpha_p\geq 0$ 
  since $\alpha_1-\alpha_1'$ is psef by assumption. 
\end{proof} 
As a consequence, we get the following useful uniformity result. 
\begin{cor}\label{bounded} 
  Let $\alpha_1,\dots,\alpha_n\in\CNS(\fX)$ be arbitrary  
  Cartier classes,  
  and suppose that for some $0\leq p\leq n$ $\alpha_i$ is nef  
  for $i\leq p$, and that  we are given a nef class  
  $\omega\in\CNS(\fX)$ such that $\omega\pm\alpha_i$  
  is nef for each $i>p$. Then we have  
  \begin{equation*} 
    |(\alpha_1\cdot\ldots\cdot\alpha_n)| 
    \leq C_n(\alpha_1\cdot\ldots\cdot\alpha_p\cdot\omega^{n-p}) 
  \end{equation*} 
  for some constant $C_n$ only depending on $n$.  
\end{cor} 
\begin{proof} 
  Write $\alpha_i$ as the difference of two nef classes 
  $\beta_i:=\alpha_i+\omega$ and $\omega$ for $i>p$.  Expanding out 
  the product of the $\alpha_i$, we see that it is enough to bound 
  terms of the form $(\alpha_1\cdot\ldots\cdot\alpha_p\cdot\beta_{i_1} 
  \cdot\ldots\cdot\beta_{i_k}\cdot\omega^{n-p-k})$.  But we also have 
  $\beta_i\leq 2\omega$ by assumption, thus the result follows from 
  Proposition~\ref{p:increase}. 
\end{proof} 
 
\subsection{Positive intersections of big classes} 
The aim of this section is to justify  
\begin{defi}\label{def:pos} 
  If $\alpha_1,\dots,\alpha_p\in \CNS(\fX)$ are 
  big classes, their positive intersection product 
  \begin{equation*} 
    \langle\alpha_1\cdot\ldots\cdot\alpha_p\rangle\in \Np(\fX) 
  \end{equation*} 
  is defined as the least upper bound of the set of classes 
  \begin{equation*} 
    (\alpha_1-D_1)\cdot\ldots\cdot(\alpha_p-D_p)\in \Np(\fX) 
  \end{equation*} 
  where 
  $D_i$ is an effective Cartier $\ku$-divisor on $\fX$ such that 
  $\alpha_i-D_i$ is nef. 
\end{defi} 
We shall see in Proposition~\ref{prop:approx} that one can in fact 
replace the $D_i$ by arbitrary psef Cartier classes. 
 
\smallskip 
To justify the definition above, we rely on two lemmas.  
Recall that a partially ordered set is \emph{directed} 
if any two elements can be dominated by a third. Dually, 
it is \emph{filtered} if any two elements dominate a third. 
\begin{lem}\label{lem:directed} 
  Let $\al\in \CNS(\fX)$ be a big Cartier class. Then the set 
  $\cD(\al)$ of effective Cartier $\ku$-divisors $D$ on $\fX$  
  such that $\al-D $ is nef is nonempty and filtered. 
\end{lem} 
\begin{lem}\label{lem:abstract} 
  Let $V$ be a Hausdorff topological vector space and $K$  
  a strict closed convex cone, with  
  associated partial order relation $\le$.  
  Then any nonempty directed subset $S\subset V$ that is  
  contained in a compact subset of $V$ admits  
  a least upper bound with respect to $\le$. 
\end{lem} 
The existence of the least upper bound in Definition~\ref{def:pos} is 
then obtained by applying Lemma~\ref{lem:abstract} to $V= \Np(\fX)$, 
$K =$ the psef cone and 
$S=\{(\alpha_1-D_1)\cdot\ldots\cdot(\alpha_p-D_p)\ \text{with}\ 
D_i\in\cD(\al_i)\}$.  The fact that $S$ is directed is a consequence 
of Lemma~\ref{lem:directed} and Proposition~\ref{p:increase}, and $S$ 
is contained in the compact set (see Proposition~\ref{compact}) of 
classes $\le \omega^p$ if $\omega$ is any Cartier class determined in 
a common determination of the $\al_i$, and sufficiently big so that 
$\omega\ge\al_i$ for all $i$. 
\begin{proof}[Proof of Lemma~\ref{lem:directed}] 
  Since $\alpha$ is big, one can find an effective  
  $\ku$-divisor $D$ on a determination $X_\pi$ of $\alpha$  
  such that $\alpha_\pi-D$ is ample, and this proves that  
  $\cD(\al)$ is nonempty. In order to prove 
  that it is filtered,  
  we will show that given two effective   $\ku$-divisors 
  $D_1, D_2$ such that $\alpha - D_i $ are nef,  
  the Weil $\ku$-divisor $D:=\min(D_1,D_2)$ on $\fX$,  
  defined coefficient-wise,  
  is in fact Cartier and that $\alpha-D$ is nef.  
  First, we can assume that $\alpha$ and the  
  $D_i$ are determined on $X$.  
  By homogeneity, we can also assume that the $D_i$ have  
  integer coefficients.  
  We then have the following characterization of $D$.  
  If we introduce the ideal sheaf 
  $\mathcal{I}:= \mathcal{O}_X(-D_1)+\mathcal{O}_X(-D_2)$,  
  then it is easy to 
  see that the incarnation $D_\pi$ of $D=\min(D_1,D_2)$  
  on a blow-up $X_\pi$ is the divisorial part of the ideal 
  $\mathcal{I}\mathcal{O}_{X_{\pi}}$.  
  This shows that $D$ is Cartier and determined in 
  $X_\pi$ if $\mathcal{I}\mathcal{O}_{X_{\pi}}$  
  is a principal ideal sheaf,  
  which is the case as soon as $\pi$ dominates  
  the normalized blow-up of   $\mathcal{I}$. 
  Since nefness is a closed condition, in order to show that  
  $\alpha-D$ is nef, we can add an arbitrarily  
  small ample class to $\alpha$  
  and reduce by homogeneity to the case where  
  $\alpha=c_1(L)$ for some line bundle $L$ on $X_\pi$  
  such that $L-D_i$ are globally generated on $X$.  
  It then follows that $\mathcal{O}_X(L)\otimes \mathcal{I}$  
  is also globally generated on $X$, and thus that 
  $\pi^{\ast}L-D$ is globally generated on $X_\pi$.  
  In particular, it is nef, and this concludes the proof  
  of the lemma. 
\end{proof} 
\begin{proof}[Proof of Lemma~\ref{lem:abstract}] 
  The assumption that $S$ is directed means that it is a net with 
  respect to the order relation $\le$.  Obviously, any accumulation 
  point of this net is a least upper bound for $S$; in particular 
  there can be at most one such accumulation point. But since this net 
  lives in a compact set, it admits an accumulation point, and 
  therefore converges towards the least upper bound of $S$. 
\end{proof} 
\begin{remark}  
  We have shown above that $\min(D_1,D_2)$ is a Cartier divisor  
  on $\fX$ for   any two effective Cartier $\ku$-divisors  
  $D_i$ on $\fX$. The relation $\min(D_1,D_2)=\min(D_1+E,D_2+E)-E$  
  then shows that the min (and thus also the max) of  
  arbitrary Cartier $\ku$-divisors on $\fX$ is Cartier.  
  All this fails, however, for arbitrary Cartier $\er$-divisors,  
  as can be seen already in the toric setting. 
\end{remark}  
\subsection{Continuity properties} 
We now proceed by describing some general properties of the positive 
intersection product. 
\begin{prop}\label{prop:cont}The positive intersection product   
  \begin{equation*} 
    (\al_1, \dots , \al_p) \mapsto 
    \langle\alpha_1\cdot\ldots\cdot\alpha_p\rangle\in \Np(\fX) 
  \end{equation*} 
  is symmetric, homogeneous of degree $1$, and super-additive in each 
  variable.   
  Moreover it varies continuously on the $p$-fold  
  product of the big cone of $\CNS(\fX)$. 
\end{prop} 
Note that in general the positive intersection product is not linear 
in each variable, as can be already seen for the map $\al \mapsto 
\langle \al \rangle$ on surfaces, using its interpretation in terms of 
the Zariski decomposition, see Section~\ref{sec:zar} below. 
\begin{proof} 
  Only the continuity statement is not clear. 
  Let $\alpha_i\in \CNS(\fX)$, $ 1 \le i \le p$, be big classes and 
  $\ep>0$.  Since $\alpha_i$ lies in the interior of the psef cone, we 
  have $\ep\alpha_i\geq\pm\gamma_i$ for every small enough Cartier 
  class $\gamma_i\in \CNS(\fX)$. It follows that 
  \begin{equation*} 
    (1-\ep)~\alpha_i\leq\alpha_i+\gamma_i\leq(1+\ep)~\alpha_i 
  \end{equation*} 
  and thus 
  \begin{equation*} 
    (1-\ep)^p~\langle\alpha_1\cdot\ldots\cdot\alpha_p\rangle \leq 
    \langle(\alpha_1+\gamma_1)\cdot\ldots\cdot(\alpha_p+\gamma_p)\rangle 
    \leq (1+\ep)^p~\langle\alpha_1\cdot\ldots\cdot\alpha_p\rangle 
  \end{equation*} 
  for  $\gamma_i$ small enough, which concludes the proof. 
\end{proof} 
We now extend the definition of the positive intersection  
product as follows. 
\begin{defi}\label{def:pos-psef} 
  If $\alpha_1,\dots,\alpha_p\in \CNS(\fX)$ are psef classes, 
  their positive intersection product 
  \begin{equation*} 
    \langle\alpha_1\cdot\ldots\cdot\alpha_p\rangle\in \Np(\fX) 
  \end{equation*} 
  is defined as the limit 
  \begin{equation*} 
    \lim_{\ep\to 0+}\langle(\alpha_1+\ep\omega)\cdot\ldots\cdot 
    (\alpha_p+\ep\omega)\rangle 
  \end{equation*} 
  where $\omega\in \CNS(\fX)$ is 
  \emph{any} big Cartier class. 
\end{defi} 
This definition makes sense, because the positive intersection 
products on the 
right-hand side decrease as $\ep$ decreases to $0$. In particular they 
lie in a compact subset of $\Np(\fX)$, so it is easy to see that the 
limit in question exists.  Furthermore, if $\omega'$ is another big 
class, then there exists $C>0$ such that $C^{-1}\omega\leq\omega'\leq 
C\omega$, and this shows that the limit is independent of the choice 
of the big class $\omega$. 
\begin{remark}  
  The extension of the positive intersection product 
  is not continuous up to the boundary of the psef cone in general, 
  see Example~\ref{ex:discon}. It is however upper semi-continuous in 
  the appropriate sense. 
\end{remark} 
An important property of the positive intersection product  
is that it coincides 
with the usual intersection product on nef classes. 
\begin{prop}  
  If $\alpha_1,\dots,\alpha_p\in \CNS(\fX)$ are nef classes, 
  then 
  \begin{equation*} 
    \langle\alpha_1\cdot\ldots\cdot\alpha_p\rangle 
    =\alpha_1\cdot\ldots\cdot\alpha_p. 
  \end{equation*} 
\end{prop} 
\begin{proof} 
  When the $\alpha_i$ are big and nef, 
  the divisor $D=0$ is allowed in Definition~\ref{def:pos}, thus 
  Proposition~\ref{p:increase} immediately yields the desired 
  equality. The general case follows, since $\omega$ can be chosen to 
  be big and nef in Definition~\ref{def:pos-psef}, in which case the 
  big classes $\alpha_i+\ep\omega$ are also nef. 
\end{proof} 
In view of this result and monotonicity, we get the following 
characterization of the positive intersection product of big classes. 
\begin{prop}\label{prop:approx}  
  If $\alpha_1,\dots,\alpha_p\in \CNS(\fX)$ are big Cartier classes, 
  their positive intersection product 
  $\langle\alpha_1\cdot\ldots\cdot\alpha_p\rangle\in \Np(\fX)$ is the 
  least upper bound of the set of all intersection products 
  \begin{equation*} 
    \beta_1\cdot\ldots\cdot\beta_p\in \Np(\fX) 
  \end{equation*} 
  with $\beta_i\in 
  \CNS(\fX)$ a nef class such that $\beta_i\leq\alpha_i$. 
\end{prop} 
\begin{remark}  
  We do not know if this result still holds for arbitrary psef classes 
  in general. It does hold when the $\alpha_i$ admit a Zariski 
  decomposition, see Section~\ref{sec:zar}. 
\end{remark} 
 
\subsection{Concavity properties}\label{sec:concave} 
The intersection products of nef classes satisfy several remarkable 
inequalities, whose proofs are based on the Hodge index theorem. We 
refer to~\cite[\S 1.6]{Lazbible} for their statements.  Thanks to 
Proposition~\ref{prop:approx}, we can transfer them to the positive 
intersection product on big classes, and by approximation to psef 
classes. This yields the following two results. 
\begin{theo}[Khovanskii-Teissier inequalities]\label{thm:KT} 
  If $\alpha_1, \cdots, \alpha_n\in \CNS(\fX)$ are psef classes,  
  then we have 
  \begin{equation*} 
    \langle\alpha_1\cdot\ldots\cdot\alpha_n\rangle\geq 
    \langle\alpha_1^n\rangle^{\frac 1 
      n}\cdots\langle\alpha_n^n\rangle^{\frac 1 n}. 
  \end{equation*} 
  More generally, we 
  have 
  \begin{equation*} 
    \langle\alpha_1\cdot\ldots\cdot\alpha_n\rangle\geq 
    \langle\alpha_1^p\cdot\alpha_{p+1}\cdot\ldots\cdot\alpha_n\rangle^{\frac1p} 
    \cdots 
    \langle\alpha_p^p\cdot\alpha_{p+1}\cdot\ldots\cdot\alpha_n\rangle^{\frac1p} 
  \end{equation*} 
  for every $1\leq p\leq n$. 
\end{theo} 
Note that this implies in particular that the sequence 
$k\mapsto\log\langle\alpha^k\cdot\beta^{n-k}\rangle$ is concave for 
any two psef classes $\alpha,\beta\in \CNS(\fX)$. 
\begin{theo}\label{thm:ineq-vol}  
  The function 
  $\alpha\mapsto \langle\alpha^n\rangle^{1/n}$ is 
  concave and homogeneous on the psef cone. More generally, given psef 
  classes $\alpha_{p+1},\cdots,\alpha_n$, the function 
  \begin{equation*} 
    \alpha\mapsto\langle\alpha^p\cdot\alpha_{p+1}\cdot \dots \cdot 
    \alpha_n\rangle^{1/p} 
  \end{equation*} 
  is homogeneous and concave on the psef cone of $\CNS(\fX)$. 
\end{theo} 
%
%
\section{The volume function}\label{sec:vol} 
 
\subsection{Continuity properties} 
Recall that the volume of a line bundle $L$ on $X$ is defined by 
\begin{equation*} 
  \vol(L)=\limsup_{k\to\infty}\frac{n!}{k^n}~h^0(X,kL). 
\end{equation*} 
As explained in the introduction, Fujita's theorem~\cite{Fujita} can 
be stated as follows. 
\begin{theo}[Fujita's theorem]  
  If  $L$ is any big line bundle on $X$, then 
  \begin{equation*} 
    \vol (L) = \langle L^n \rangle. 
  \end{equation*} 
\end{theo} 
Since the function  
$\alpha\mapsto\langle \alpha^n \rangle$  
is homogeneous of degree $n$ and continuous on the big cone,  
it follows that it coincides with the (necessarily unique)  
continuous and $n$-homogeneous extension of the volume function  
to the big cone as constructed (using different arguments) 
in~\cite[\S2.2.C]{Lazbible}.  Now positive intersection  
products are not continuous up to the boundary of the psef cone  
in general, see Example~\ref{ex:discon},  
and it is thus remarkable that the following holds: 
\begin{theo}\label{bound-cont} 
  The function $\alpha\mapsto\langle\alpha^n\rangle$  
  is strongly continuous on the psef cone of $\CNS(\fX)$ 
  and vanishes on its boundary (and only there).  
\end{theo} 
This result follows from~\cite[Cor~2.2.45]{Lazbible}, and is also 
proved for arbitrary $(1,1)$-classes using analytic techniques 
in~\cite{Boucksom02}. Note that since 
$\alpha\mapsto\langle\alpha^n\rangle$ is nonnegative and upper 
semi-continuous up to the boundary of the psef cone, its continuity 
follows in fact from its vanishing on the boundary. All in all we get 
that the function $\vol:\CNS(\fX)\to\er$ defined by 
\begin{equation*} 
  \vol(\alpha):=\langle\alpha^n\rangle 
\end{equation*} 
when $\alpha$ is psef, and  
\begin{equation*} 
  \vol(\alpha):=0 
\end{equation*} 
when $\alpha$ is not psef, is continuous and coincides with the volume 
function defined in~\cite[\S2.2.C]{Lazbible}. 
 
\subsection{Proof of Theorem A} 
Let us first recall the following fundamental Morse-type inequality: 
\begin{prop}\label{morse} For any two nef Cartier classes 
  $A,B\in\CNS(\fX)$, we have 
  \begin{equation*} 
    \vol(A-B)\geq(A^n)-n(A^{n-1}\cdot B) 
  \end{equation*} 
\end{prop} 
We refer to~\cite[Proof of Theorem~2.2.15]{Lazbible} for an elementary 
algebraic proof using the interpretation of the volume in terms of 
growth of sections. As an immediate consequence, we get the following 
sub-linear control of the volume near a nef class: 
\begin{cor}\label{sublin} 
  Let $\beta\in\CNS(\fX)$ be a nef class, and $\gamma\in\CNS(\fX)$  
  an arbitrary Cartier class. If $\omega\in\CNS(\fX)$ is a given nef 
  and big class such that $\beta\leq\omega$ and $\omega\pm\gamma$  
  is nef, then we have  
  \begin{equation*} 
    \vol(\beta+t\gamma)\geq(\beta^n)+nt(\beta^{n-1}\cdot\gamma)-Ct^2 
  \end{equation*} 
  for every $0\leq t\leq 1$ and some constant $C>0$ only  
  depending on  $(\omega^n)$. 
\end{cor}  
\begin{proof}  
  We claim  
  that $(\beta+t\gamma)^n=(\beta^n)+nt(\beta^{n-1}\cdot\gamma)+O(t^2)$ 
  for $0\leq t\leq 1$, with $O$ only depending on $(\omega^n)$. 
  Indeed, $O$ is controlled by $(\beta^k\cdot\gamma^{n-k})$, 
  $k=0,\dots,n-2$, 
  and thus by $(\beta^k\cdot\omega^{n-k})$ thanks to 
  Proposition~\ref{bounded}, which in turn is bounded by $(\omega^n)$ 
  according to Proposition~\ref{p:increase} since 
  $\beta\leq\omega$. Now we write $\beta+t\gamma$ as the difference of 
  the two nef classes $A:=\beta+t(\gamma+\omega)$ and 
  $B:=t\omega$. Then we also have 
  \begin{equation*} 
    (\beta+t\gamma)^n=(A-B)^n=(A^n)-n(A^{n-1}\cdot B)+O(t^2) 
  \end{equation*} 
  with $O$ only depending on $(\omega^n)$. Indeed, this $O$ is 
  controlled by $(A^k\cdot\omega^{n-k})$, $k=0,\dots,n-2$, and we have 
  $A\leq 3\omega$, so we again get a control in terms of $(\omega^n)$ 
  only.  All in all we thus have 
  \begin{equation*} 
    (A^n)-n(A^{n-1}\cdot B) 
    =(\beta^n)+nt(\beta^{n-1}\cdot\gamma)+O(t^2) 
  \end{equation*} 
  and the result now follows from an application of 
  Proposition~\ref{morse} to $\beta+t\gamma=A-B$. 
\end{proof} 
We are now in a position to prove our main result, Theorem A. We thus 
consider a big Cartier class $\alpha\in\CNS(\fX)$ and an arbitrary 
Cartier class $\gamma\in\CNS(\fX)$, and fix a sufficiently nef and big 
class $\omega\in\CNS(\fX)$ such that $\alpha\leq\omega$ and 
$\omega\pm\gamma$ is nef. If $\beta\leq\alpha$ is a nef Cartier class, 
then \emph{a fortiori} $\beta\leq\omega$, and we deduce from 
Corollary~\ref{sublin} that 
\begin{equation*} 
  \vol(\alpha+t\gamma)\geq 
  \vol(\beta+t\gamma)\geq(\beta^n)+nt(\beta^{n-1}\cdot\gamma)-Ct^2 
\end{equation*} 
for every $0\leq t\leq 1$ and some constant $C>0$ only depending on 
$(\omega^n)$. Taking the supremum over all such nef classes 
$\beta\leq\alpha$ yields 
\begin{equation*} 
  \vol(\alpha+t\gamma) 
  \geq\vol(\alpha)+nt\langle\alpha^{n-1}\rangle\cdot\gamma-Ct^2, 
\end{equation*} 
for some $C>0$ only depending on $(\omega^n)$. This holds for every 
$0\leq t\leq 1$, and in fact also for every $-1\leq t\leq 1$, merely 
by replacing $\gamma$ by $-\gamma$. Exchanging the roles of 
$\alpha+t\gamma\leq 2\omega$ and $\alpha=(\alpha+t\gamma)+t(-\gamma)$, 
this yields 
\begin{equation*} 
  \vol(\alpha)\geq\vol(\alpha+t\gamma)-nt\langle(\alpha+t\gamma)^{n-1}\rangle 
  \cdot\gamma-Ct^2, 
\end{equation*} 
for some possibly larger $C>0$ also depending only on 
$(\omega^n)$. The combination of these two inequalities immediately 
shows that 
\begin{equation*} 
  \left.\frac{d}{dt}\right|_{t=0}\vol(\alpha+t\gamma) 
  =n\langle\alpha^{n-1}\rangle\cdot\gamma 
\end{equation*} 
as desired, because $\langle(\alpha+t\gamma)^{n-1}\rangle$ converges 
to $\langle\alpha^{n-1}\rangle$ by Proposition~\ref{prop:cont}, 
$\alpha$ being big. We have thus shown that the volume admits a 
directional derivative in any direction, and that this derivative is 
given by the linear form $n\,\langle\alpha^{n-1}\rangle\in 
\NS(X)^{\ast}$. On the big cone, this form varies continuously, so 
$\vol$ is of class $\mathcal{C}^1$. 
\begin{remark} 
  In fact, we have proved that the volume function admits a 
  directional derivative in any direction in the infinite dimensional 
  space $\CNS(\fX)$, which is induced by a continuous linear form on 
  $\CNS(\fX)$. In particular, the restriction of $\vol$ to any finite 
  dimensional subspace is $\mathcal{C}^1$. 
\end{remark} 
 
As a consequence of Theorem A, we get the following orthogonality 
property, which was the key point in the characterization of 
pseudo-effectivity in~\cite{BDPP}. 
\begin{cor}\label{ortho}  
  For any psef class $\alpha\in\CNS(\fX)$, we have 
  \begin{equation*} 
    \langle\alpha^n\rangle=\langle\alpha^{n-1}\rangle\cdot\alpha. 
  \end{equation*} 
\end{cor} 
\begin{proof} 
  It is enough to show this  
  when $\alpha$ is big. Applying Theorem A to $\gamma:=\alpha$ yields 
  that $n\alpha\cdot\langle\alpha^{n-1}\rangle$ coincides with the 
  derivative at $t=0$ of 
  $\langle(\alpha+t\alpha)^n\rangle=(1+t)^n\langle\alpha^n\rangle$, 
  which is of course nothing but $n\langle\alpha^n\rangle$. 
\end{proof} 
 
 
\subsection{The Diskant inequality} 
In this section, we prove the Diskant inequality (Theorem~F), and its 
consequences Theorem~D and Corollary~E on the characterization of the 
equality case in the Khovanskii-Teissier inequalities. 
 
Before starting the proof, we formalize one ingredient in the 
statement of Diskant inequality,  
and define the \emph{slope} of a big class 
with respect to another.  Recall that since the big cone is the 
interior of the psef cone by Proposition~\ref{prop:recall}, if $\al$ 
and $\beta$ are big classes, there exists $t>0$ such that 
$t\beta\le\alpha$. 
 
\begin{defi} 
  The \emph{slope} of $\beta$ with respect to $\alpha$ is 
  defined as 
  \begin{equation*} 
    s=s(\alpha,\beta)=\sup\{t>0\ |\ \alpha\ge t\beta\}. 
  \end{equation*} 
\end{defi} 
Since the psef cone is closed, we have 
$\alpha\geq s\beta$, and $\alpha-t\beta$ is big for $t<s$. 
Note that $\alpha=\beta$ iff $s(\alpha,\beta)=s(\beta,\alpha)=1$. 
\begin{proof}[Proof of Theorem~D] 
  Since $k\mapsto\log(\alpha^k\cdot\beta^{n-k})$ is concave,~(i) 
  and~(ii) are equivalent. Moreover,~(iii) trivially implies~(i), so 
  it only remains to prove that~(i) implies~(iii). By homogeneity we 
  may assume that $(\alpha^n)=(\beta^n)=1$, which implies 
  $s(\alpha,\beta)\le1$ by Proposition~\ref{p:increase}. From~(i) we 
  get $(\alpha^{n-1}\cdot\beta)=1$ so Diskant's inequality gives 
  $s(\alpha,\beta)=1$. By symmetry we get $s(\beta,\alpha)=1$. Hence 
  $\alpha=\beta$. 
\end{proof} 
\begin{proof}[Proof of Corollary~E] 
  Pick $\alpha,\beta$ big and nef with 
  $((\alpha+\beta)^n)^{\frac1n}=(\alpha^n)^{\frac1n}+(\beta^n)^{\frac1n}$. 
  Then 
  \begin{align*} 
    ((\alpha+\beta)^n) 
    &=\sum_{k=0}^n\binom{n}{k}(\alpha^k\cdot\beta^{n-k})\\ 
    &\le\sum_{k=0}^n\binom{n}{k} 
    (\alpha^n)^{\frac{k}{n}}(\beta^n)^{\frac{n-k}{n}} 
    =((\alpha^n)^{\frac1n}+(\beta^n)^{\frac1n})^n, 
  \end{align*} 
  so the function $k\mapsto\log(\alpha^k\cdot\beta^{n-k})$  
  must be affine. 
  By Theorem~D, this implies that $\alpha$ and $\beta$ 
  are proportional. 
\end{proof} 
\begin{proof}[Proof of Theorem~F] 
  Set $\alpha_t:=\alpha-t\beta$ for $t\geq 0$. 
  By the definition of the slope $s=s(\alpha,\beta)$, 
  $\alpha_t$ is big iff $t<s$. 
  By Theorem~A, 
  $f(t):=\vol(\alpha_t)$ is differentiable for $t<s$, with 
  $f'(t)=-n\langle\alpha_t^{n-1}\rangle\cdot\beta$. 
  We also have $f(0)=(\alpha^n)$ and 
  $f(t)\to 0$ as $t\to s_-$ by continuity 
  so it follows that 
  \begin{equation*} 
    (\alpha^n)=n\int_{t=0}^s\langle\alpha_t^{n-1}\rangle\cdot\beta\,dt. 
  \end{equation*} 
  If $\gamma_t\in \CNS(\fX)$ is a nef class with 
  $\gamma_t\leq\alpha_t=\alpha-t\beta$, then 
  $\gamma_t+t\beta\leq\alpha$ implies 
  \begin{equation*} 
    (\gamma_t^{n-1}\cdot\beta)^{\frac1{n-1}}+t(\beta^n)^{\frac1{n-1}} 
    \leq(\alpha^{n-1}\cdot\beta)^{\frac1{n-1}} 
  \end{equation*} 
  by Theorem~\ref{thm:ineq-vol}. 
  Taking the supremum over all such nef classes $\gamma_t$ yields 
  \begin{equation*} 
    (\langle\alpha_t^{n-1}\rangle\cdot\beta)^{\frac1{n-1}} 
    +t(\beta^n)^{\frac1{n-1}} 
    \leq(\alpha^{n-1}\cdot\beta)^{\frac1{n-1}}, 
  \end{equation*} 
  by Proposition~\ref{prop:approx}. 
  We thus obtain 
  \begin{equation*} 
    (\alpha^n)\leq n\int_{t=0}^s\left((\alpha^{n-1}\cdot\beta)^{1/n-1} 
      -t(\beta^n)^{1/n-1}\right)^{n-1}dt 
  \end{equation*} 
  and the result follows since 
  \begin{equation*} 
    \frac{d}{dt}\left((\alpha^{n-1}\cdot\beta)^{\frac1{n-1}} 
      -t(\beta^n)^{\frac1{n-1}}\right)^n\!  
    =n(\beta^n)^{\frac1{n-1}} 
    \left((\alpha^{n-1}\cdot\beta)^{\frac1{n-1}} 
      -t(\beta^n)^{\frac1{n-1}}\right)^{n-1} 
  \end{equation*} 
\end{proof} 
 
 
\subsection{Positive intersection and Zariski decomposition}\label{sec:zar} 
Although the discussion to follow is strictly speaking not necessary  
for the understanding of the rest of the article,  
we would like to emphasize at this point that there is a  
very close relationship between positive intersection products and  
Zariski-type decompositions of psef classes.   
A survey of the different definitions of Zariski decompositions  
that have been proposed in higher dimension (and which all coincide  
for big classes) is given in~\cite{prok}, and a very complete account 
on the existence problem can be found in~\cite{nakayama}. 
 
For a psef class $\alpha \in \CNS(\fX)$, set  
$P(\alpha):=\langle\alpha\rangle\in \NS(\fX)$, so that  
$P(\alpha)\leq\alpha$. By definition, any nef Cartier class  
$\beta\in\CNS(\fX)$ such that $\beta\leq\alpha$ already satisfies  
$\beta\leq P(\alpha)$. One easily deduces  
that the incarnation of the Weil class $P(\alpha)$ on each  
$X_\pi$ coincides with the positive part of $\alpha_\pi$  
in its divisorial Zariski decomposition, as first introduced by  
Nakayama (see~\cite{nakayama}, where it is called $\sigma$-decomposition),  
and independently for arbitrary $(1,1)$-classes via analytic tools  
by the first author in~\cite{bouck}. This means that the collection  
of all the divisorial Zariski decompositions in all models  
gives rise to a decomposition  
$\alpha =P(\alpha) + N(\alpha)$ in $\NS(\fX)$,  
where $N(\alpha)$ is (the class of) an effective Weil divisor on $\fX$.  
 
We will say that $\alpha$ \emph{admits a Zariski decomposition}  
if $P(\alpha)$ is a Cartier class on $\fX$.  
This is equivalent to requiring that the  
positive part of $\alpha_\pi$ in its divisorial Zariski decomposition  
on $X_\pi$ is nef for some $\pi$, i.e.\ that $P(\alpha)$ is a nef  
Cartier class on $\fX$, a situation called 
"generalized Fujita decomposition" in~\cite{prok}.  
When psef classes $\alpha_1,\dots,\alpha_p$ admit a Zariski decomposition,  
it is clear that 
\begin{equation*} 
  \langle\alpha_1\cdot\ldots\cdot\alpha_p\rangle 
  =P(\alpha_1)\cdot\ldots\cdot P(\alpha_p) 
\end{equation*} 
The Fujita theorem and the orthogonality relation  
(Corollary~\ref{ortho}) can thus be formulated respectively as  
\begin{equation*} 
  \vol (\alpha) =P(\alpha)^n\ \text{and}\ P(\alpha)^{n-1} \cdot 
  N(\alpha)=0 
\end{equation*} 
when $\alpha$ admits a Zariski decomposition (the terminology 
"orthogonality relation" in fact stems from this second equality). 
 
In dimension $n=2$, it follows from classical results that  
Zariski decompositions always exist.  
However, in higher dimensions, psef (or even big) Cartier classes do 
not admit a Zariski decomposition in general.  A counter-example is 
constructed in~\cite{nakayama}, in which $\alpha$ is the big Cartier 
class induced by the tautological line bundle on 
$X:=\mathbf{P}_S(L_1\oplus L_2\oplus L_3)$, the $L_i$ being 
appropriate line bundles on an abelian surface $S$.

In the general case where classes do not necessarily admit Zariski 
decompositions, one can interpret the construction of positive 
intersection products as making sense of the intersection 
$P_1\cdot\ldots\cdot P_p$ of the Weil classes $P_i=P(\alpha_i)$, even 
though it is definitely not possible to make sense of the intersection 
of arbitrary Weil classes on $\fX$. 
 
\medskip 
Using the existence and characterization of Zariski decompositions on 
surfaces recalled above, we can give the following counter-example to 
continuity of positive intersection products up to the boundary of the 
psef cone. 
\begin{example}\label{ex:discon} 
  Let $X$ be any projective surface with infinitely many  
  exceptional curves,   i.e.\ irreducible curves $C_k$ such that  
  $(C_k^2)<0$ (for instance the blow-up of  
  $\mathbf P^2$ at $9$ points).  
  If $\omega$ is a given ample class on $X$ and 
  $t_k:=(\omega\cdot C_k)^{-1}>0$, then $\beta_k:=t_kC_k$  
  is bounded in 
  $\NS(X)$, thus we can assume that it converges to a non-zero class 
  $\beta\in \NS(X)$.  Since the $C_k$ are distinct, the limit 
  $\beta$ is nef, and thus $P(\beta)=\beta$.  
  But as $C_k$ is contractible, we   get $P(\beta_k)=0$,  
  which shows that $\alpha\mapsto P(\alpha)$ is not  
  continuous at $\beta$. If $\alpha$ is any ample class on $X$,  
  we have $\langle\alpha\cdot\beta_k\rangle=\alpha\cdot P(\beta_k)=0$,  
  whereas $\langle\alpha\cdot\beta\rangle=\alpha\cdot\beta\neq 0$,  
  and this shows indeed that the positive intersection product  
  is not continuous at $(\alpha,\beta)$.  
  Note that $\beta$ lies on the boundary of the psef cone.  
\end{example} 
 
Beyond surfaces, toric varieties constitute an important class of 
varieties on which Zariski decompositions always exist.  Indeed, if 
$X$ is a projective toric variety, then the nef part of a toric psef 
class $\alpha\in\CNS(\fX_{\mathrm{tor}})$ is just the toric nef class 
$P(\alpha)\in\CNS(\fX_{\mathrm{tor}})$ associated to the Newton 
polytope $\mathrm{Nw}(\alpha)$ of $\alpha$ 
(cf.~Section~\ref{sec:toric}).  In other words, if $g$ is the 
homogeneous function on $\er^n$ corresponding to $\alpha$, which is 
piecewise linear with respect to some rational fan decomposition 
$\Sigma$ of $\er^n$, then its nef part corresponds to its convex 
minorant, i.e.\ the largest homogeneous and convex function $h\leq g$. 
This function $h$ is also piecewise linear with respect to some 
refinement $\Sigma'$ of the fan $\Sigma$. 
This means that if $\alpha$ is determined on some toric blow-up  
$X_\pi$ of $X$, the Cartier class $P(\alpha)$ is determined on  
some higher toric blow-up $\pi'\geq\pi$, which cannot be taken  
to be $\pi$ in general (as was the case for surfaces).  
On the other hand, for each fixed $\pi$ the map $\alpha\mapsto P(\alpha)$  
is piecewise linear on the psef cone of $\NS(X_{\pi})$ with respect  
to the Gel'fand-Kapranov-Zelevinskij decomposition of~\cite{oda}.  
This implies in particular that the volume $\vol(\alpha)=P(\alpha)^n$  
is piecewise polynomial on the psef cone of $\NS(X_\pi)$,  
as explained in~\cite{ELMNP2}. 
 
We conclude this section by describing the relationship between  
positive products and Brunn-Minkowski theory.  
If $L$ is a big line bundle on some toric blow-up $X_\pi$, then  
$h^0(X,L)$ is equal to the number of integral points in the  
Newton polytope $\mathrm{Nw}(L) \subset \er^n$. 
Since $\mathrm{Nw}(kL)=k\mathrm{Nw}(L)$,  
it follows that $\frac{\vol(L)}{n!}=\lim_{k\to\infty}h^0(X,kL)/k^n$  
equals the Euclidean volume of $\mathrm{Nw}(L)$.   
This then extends by linearity and continuity to show that if  
$\alpha_i\in \CNS(\fX_{\mathrm{tor}})$, $i=1,\dots,n$  
are psef toric classes, the positive intersection product 
$\langle\alpha_1\cdot\ldots\cdot\alpha_n\rangle 
=P(\alpha_1)\cdot\ldots\cdot P(\alpha_n)$ 
is nothing but the mixed volume of their  
Newton polytopes $\mathrm{Nw}(\alpha_i)$ as defined in  
Brunn-Minkowski theory (up to the factor $n!$),  
see~\cite{Schneider}.

 
\section{Restricted volumes}\label{sec:res} 
 
 
\subsection{Restriction of positive intersection products to divisors.} 
>From now on, we assume that $X$ is a projective \emph{normal} variety, 
and let $D$ be a prime divisor on $X$. In this situation, we want to 
define a restriction map from classes on the Riemann-Zariski space 
$\fX$ of $X$ to classes on the Riemann-Zariski space $\fD$ of $D$. For 
Cartier classes, this is easily done as follows. Since $X$ is normal, 
any blow-up of $X$ is an isomorphism over the generic point of $D$, 
and can be dominated by a blow-up $\pi:X_\pi\to X$ such that the 
strict transform $D_\pi$ of $D$ on $X_\pi$ is smooth, by embedded 
resolution of singularities. The corresponding system of restriction 
maps $N^p(X_\pi)\to N^p(D_\pi)$ is compatible under pull-back, and 
thus defines a continuous restriction map on Cartier classes 
\begin{equation*} 
  \CNp(\fX)\to\CNp(\fD) 
\end{equation*} 
\begin{equation*} 
  \alpha\mapsto\alpha|_{\fD}. 
\end{equation*} 
It is however not possible to extend this map to a continuous linear 
map on Weil classes $N^p(\fX)\to N^p(\fD)$ 
$\alpha\mapsto\alpha|_{\fD}$ in general. Indeed, such a continuous 
extension is necessarily unique by density of Cartier classes, and 
writing a given Weil class $\alpha\in N^p(\fX)$ as a limit of Cartier 
classes $\alpha=\lim_\pi \alpha_\pi$ as in Lemma~\ref{dense} shows 
that $\alpha\mapsto\alpha|_{\fD}$ should satisfy 
\begin{equation*} 
  (\alpha|_{\fD})_{\pi}=\alpha_{\pi}|_{D_{\pi}} 
\end{equation*} 
for every blow-up $\pi$ of $X$. But this relation already fails for 
Cartier classes, as can be seen for instance if $X$ is $\CP^2$, $D$ is 
a line and $\alpha$ is the Cartier class on $\fX$ determined by the 
strict transform of $D$ on the blow-up $X_\pi$ of $X$ at a point of 
$D$. 
 
The goal of this section will be to show that it is possible to 
define the restriction to $\fD$ of positive intersection products 
$\langle\alpha_1\cdot\ldots\cdot\alpha_p\rangle\in N^p(\fX)$ of psef 
classes $\alpha_i\in\CNS(\fX)$ under a suitable positivity assumption 
with respect to $D$. The main point is that the positive intersection 
product of big classes $\alpha_i\in\CNS(\fX)$ is by definition a 
monotone limit of Cartier classes of the form 
$\beta_1\cdot\ldots\cdot\beta_p$ with $\beta_i\leq\alpha_i$ a nef 
Cartier class.  We can therefore try to define the restriction of 
$\langle\alpha_1\cdot\ldots\cdot\alpha_p\rangle$ to $\fD$ as the 
monotone limit of the Cartier classes 
$(\beta_1\cdot\ldots\cdot\beta_p)|_{\fD}$. The trouble is that 
monotonicity is meant with respect to pseudo-effectivity, which is in 
general destroyed upon restricting to $\fD$. We thus introduce the 
following definitions. 
\begin{defi}  
  If $D$ is a prime divisor on a normal projective variety 
  $Y$, we say that a class $\alpha\in \NS(Y)$ is \emph{$D$-psef} and 
  write $\alpha\ge_D 0$ iff it belongs to the closed convex cone 
  generated by effective divisors whose support does not contain $D$. 
  A class $\alpha$ is \emph{$D$-big} if it lies in the interior  
  of the $D$-psef cone. 
\end{defi} 
It is clear that if $\alpha$ is $D$-psef then $\alpha|_D$ is 
psef. Note also that a $D$-psef class is in particular psef, and thus 
that a $D$-big class is in particular big. 
\begin{remark} 
  If $\alpha\in\NS(Y)$ is psef, then:  
  \begin{itemize} 
  \item 
    $\alpha$ is $D$-psef iff $D$ is not contained in the restricted 
    base locus (or non-nef locus) of $\alpha$, 
    see~\cite{ELMNP2,bouck}; 
  \item 
    $\alpha$ is $D$-big iff $D$ is not contained in the augmented base 
    locus (or non-ample locus) of $\alpha$, see~\cite{ELMNP2,bouck}. 
  \end{itemize} 
\end{remark} 
If $\mu:Y'\to Y$ is a birational morphism, and if $D'$ denotes the 
strict transform of $D$, then $\mu^{\ast}\alpha$ is $D'$-psef iff 
$\alpha\in \NS(Y)$ is $D$-psef. Considering again a normal projective 
variety $X$ and a prime divisor $D$ on $X$, it follows that the 
following definition makes sense: 
\begin{defi}  
  A class 
  $\alpha\in \CNS(\fX)$ is $\fD$-psef (resp. $\fD$-big) iff there 
  exists a determination $\pi$ of $\alpha$ such that $\alpha_\pi\in 
  \NS(X_\pi)$ is $D_\pi$-psef (resp. $D_\pi$-big). 
\end{defi} 
We will write $\alpha\geq_{\fD}0$ if $\alpha\in\CNS(\fX)$ is 
$\fD$-psef. Note that a class in $\CNS(\fX)$ is $\fD$-big iff it 
belongs to the interior of the $\fD$-psef cone in the strong 
topology. 
 
\smallskip 
Proceeding as in Section~\ref{sec:pos}, one shows that the following definition 
makes sense. 
\begin{defi}\label{def:above} 
  If $\alpha_1,\dots,\alpha_p\in \CNS(\fX)$ are $\fD$-big 
  classes, one defines their restricted  
  positive intersection product on $\fD$ 
  \begin{equation*} 
    \langle\alpha_1\cdot\ldots\cdot\alpha_p\rangle|_{\fD}\in \Np(\fD) 
  \end{equation*} 
  as the least upper bound of the set of classes 
  \begin{equation*} 
    (\beta_1\cdot\ldots\cdot\beta_p)|_{\fD}\in \Np(\fD) 
  \end{equation*} 
  where 
  $\beta_i\in \CNS(\fX)$ is a nef class such that 
  $\beta_i\leq_\fD\alpha_i$. 
\end{defi} 
\begin{remark}  
  The $\fD$-big classes $\alpha_i\in \CNS(\fX)$ restrict to big 
  classes $\alpha_i|_{\fD}\in \CNS(\fD)$, so we can also consider 
  their positive intersection on $\fD$, to wit 
  $\langle\alpha_1|_{\fD}\cdot\ldots\cdot\alpha_p|_{\fD} \rangle\in 
  \Np(\fD)$. It is the least upper bound of the set of classes 
  $\gamma_1\cdot\ldots\cdot\gamma_p\in \Np(\fD)$ where $\gamma_i\in 
  \CNS(\fD)$ is nef and such that $\gamma_i\leq \alpha_i|_{\fD}$.  The 
  point of Definition~\ref{def:above} above is that we only consider 
  nef classes $\gamma_i\leq \alpha_i|_{\fD}$ that are 
  \emph{restrictions} to $\fD$ of nef classes 
  $\beta_i\leq_\fD\alpha_i$ on $\fX$. In particular, it follows that 
  \begin{equation*} 
    \langle\alpha_1\cdot\ldots\cdot\alpha_p\rangle|_{\fD} 
    \leq\langle \alpha_1|_{\fD}\cdot\ldots\cdot\alpha_p|_{\fD}\rangle 
  \end{equation*} 
  in $\Np(\fD)$, but equality does not hold in general. 
\end{remark} 
Since the restricted positive intersection product is homogeneous and 
increasing in each variable (with respect to $\geq_\fD$), continuity 
holds as in Proposition~\ref{prop:cont}. 
\begin{prop}\label{cont-rest} 
  The restricted positive intersection product 
  $\langle\alpha_1\cdot\ldots\cdot\alpha_p\rangle|_{\fD}\in \Np(\fD)$ 
  depends continuously on the $\fD$-big classes $\alpha_i\in 
  \CNS(\fX)$. 
\end{prop} 
Again the definition can be extended to $\fD$-psef classes 
$\alpha_i\in \CNS(\fX)$ by setting 
\begin{equation*} 
  \langle\alpha_1\cdot\ldots\cdot\alpha_p\rangle|_{\fD} 
  :=\lim_{\ep\to0+} 
  \langle( 
  \alpha_1+\ep\omega)\cdot\ldots\cdot(\alpha_p+\ep\omega) 
  \rangle|_{\fD} 
\end{equation*} 
for $\omega\in \CNS(\fX)$ a $\fD$-big class.  Indeed, the limit in 
question does not depend on the choice of $\omega$. It depends upper 
semi-continuously on the $\fD$-psef classes $\alpha_i$, but continuity 
does not hold up to the boundary in general. 
\begin{prop}\label{prop:car-rest} 
  If $\alpha_1,\dots,\alpha_p\in \CNS(\fX)$ are nef, then 
  \begin{equation*} 
    \langle\alpha_1\cdot\ldots\cdot\alpha_p\rangle|_{\fD}= 
    (\alpha_1\cdot\ldots\cdot\alpha_p)|_{\fD} 
  \end{equation*} 
\end{prop} 
We also note that the concavity properties of 
Section~\ref{sec:concave} also hold in this setting, again because the 
Khovanskii-Teissier inequalities for nef classes hold on $\fD$. 
 
\smallskip 
We will use the following easy inequality:  
\begin{prop}\label{compar}  
  Let $D$ be a prime divisor on $X$ and assume that  
  $D$ is also Cartier. If $\alpha_1,\dots,\alpha_{n-1}\in\CNS(\fX)$  
  are $\fD$-psef classes, then we have 
  \begin{equation*} 
    \langle\alpha_1\cdot\ldots\cdot\alpha_{n-1}\rangle|_{\fD} 
    \leq\langle\alpha_1\cdot\ldots\cdot\alpha_{n-1}\rangle\cdot D. 
  \end{equation*} 
\end{prop} 
\begin{proof} 
  It is enough consider the case of $\fD$-big classes $\alpha_i$. 
  Pick $\beta_i\in \CNS(\fX)$  a nef class with 
  $\beta_i\leq_\fD\alpha_i$.  
  If $\pi$ is a determination of $\beta_i$, then 
  \begin{multline*} 
    (\beta_1\cdot\ldots\cdot\beta_{n-1})|_{\fD} 
    =\beta_{1,\pi}\cdot\ldots\cdot\beta_{n-1,\pi}\cdot D_\pi\leq\\ 
    \leq \beta_{1,\pi}\cdot\ldots\cdot\beta_{n-1,\pi}\cdot\pi^{\ast}D 
    \leq\langle\alpha_1\cdot\ldots\cdot\alpha_{n-1}\rangle\cdot D, 
  \end{multline*} 
  where $D_\pi$ denotes as before the 
  strict transform of $D$ on $X_\pi$.  
  The desired inequality follows by taking   the supremum  
  over all such nef classes $\beta_i$. 
\end{proof} 
We can now give the following characterization of $\fD$-big classes. 
\begin{theo}\label{charac}  
  If $\alpha\in\CNS(\fX)$ is big and $\fD$-psef, then  
  $\alpha$ is $\fD$-big iff $\langle\alpha^{n-1}\rangle|_{\fD}>0$. 
\end{theo} 
\begin{proof}  
  If $\alpha$ is $\fD$-big, then $\alpha\geq_{\fD}\omega$ is for  
  some nef class $\omega\in\CNS(\fX)$ which is  
  determined by an ample class on $X$. We thus have  
  $\langle\alpha^{n-1}\rangle|_{\fD}\geq(\omega|_D)^{n-1}>0$.  
  Conversely, assume that $\alpha$ is big and $\fD$-psef  
  but not $\fD$-big, and let us show that 
  $\langle\alpha^{n-1}\rangle|_{\fD}=0$.  
  Upon replacing $X$ by some higher birational model,  
  we can also assume that $\alpha$ is determined on $X$,  
  and that $X$ is smooth, so that $D$ is in particular Cartier.  
  In view of the preceding proposition, the result will thus follow  
  from the following lemma.  
\end{proof} 
\begin{lem}\label{orthobis}  
  If $D$ is Cartier on $X$ and $\alpha$ is a big class determined  
  on $X$ which is not $D$-big,  
  then $\langle\alpha^{n-1}\rangle\cdot D=0$. 
\end{lem} 
\begin{proof} 
  If $\omega$ is an ample class on $X$, then $\alpha-\ep\omega$  
  is big but not $D$-psef for every small enough $\ep>0$.  
  By continuity, it is thus enough to show that 
  $\langle\alpha^{n-1}\rangle\cdot D=0$ for any big class  
  $\alpha$ on $X$ which is not $D$-psef.  
  But in that case, $D$ is contained in the non-nef  
  locus of $\alpha$, i.e.\ $\alpha\geq P+tD$ for some $t>0$  
  if we denote by $P$ the incarnation on $X$  
  of $\langle\alpha\rangle$ (which coincides with  
  the positive part in the divisorial Zariski decomposition  
  of $\alpha$, cf.~\cite{bouck}). It follows that 
  $\langle\alpha^{n-1}\rangle\cdot\alpha 
  \geq\langle P^{n-1}\rangle\cdot P 
  +\langle\alpha^{n-1}\rangle\cdot tD$. 
  But we also have  
  $\langle P^{n-1}\rangle\cdot P=\langle P^n\rangle 
  =\langle\alpha^n\rangle=\langle\alpha^{n-1}\rangle\cdot\alpha$ 
  by the orthogonality relation (Corollary~\ref{ortho})  
  applied to $P$ and $\alpha$, and thus we get  
  $\langle\alpha^{n-1}\rangle\cdot D=0$ as claimed. 
\end{proof} 
As noticed after Theorem~\ref{bound-cont}, this implies: 
\begin{cor}  
  The function $\alpha\mapsto\langle\alpha^{n-1}\rangle|_{\fD}$  
  is continuous on the cone of big and $\fD$-psef classes of 
  $\CNS(\fX)$. 
\end{cor} 
 
\begin{remark}  
  The theorem fails when $\alpha$ is not big. For instance,  
  if $X$ is a smooth   surface, $D$ is ample and  
  $\alpha\in\CNS(\fX)$ is nef but not big, then 
  $\langle\alpha\rangle|_{\fD}=\alpha\cdot D>0$  
  as soon as $\alpha$ is non-zero, but $\alpha$ is not  
  $\fD$-big since it is not even big.   
\end{remark} 
\begin{remark}  
  The lemma fails if $\alpha$ is not determined on $X$.  
  For instance, if $X$ is a ruled surface, $D$ is a ruling  
  and $\pi:X_{\pi}\to X$ is the blow-up of $X$ at a point $p\in D$,  
  with exceptional divisor $E$, then the strict transform  
  $D_{\pi}$ of $D$ can be blown-down by $\mu:X_{\pi}\to X'$.  
  If $\alpha_{\pi}:=\mu^{\ast}\omega$ for some ample  
  class $\omega$ on $X'$, then clearly the Cartier class  
  $\alpha\in\CNS(\fX)$ determined by $\alpha_{\pi}$  
  is nef and big, but not $\fD$-big. However, we have  
  \begin{equation*} 
    \alpha\cdot D=\alpha_{\pi}\cdot(D_\pi+E)=\omega\cdot E'>0 
  \end{equation*} 
  with $E'$ the image of $E$ on $X'$.   
\end{remark} 
 
\subsection{Restricted volumes.} 
Let again $X$ be a normal projective variety and $D$ be a prime 
divisor. Given a line bundle $L$ on $X$, we will denote by 
$h^0(X|D,L)$ the rank of the restriction map 
\begin{equation*} 
  H^0(X,L)\to H^0(D,L|_{D}). 
\end{equation*} 
Recall from the introduction that the 
restricted volume of $L$ on $D$ is then defined as 
\begin{equation*} 
  \vol_{X|D}(L):=\limsup_{k\to\infty}\frac{(n-1)!}{k^{n-1}}~h^0(X|D,kL). 
\end{equation*} 
It is thus the growth coefficient of the number of sections of 
$\oh_D(kL)$ on $D$ that extend to $X$.  
As stated in the introduction, we have 
\begin{theo}[Generalized Fujita's Theorem,~\cite{ELMNP3,Takayama}] 
  \label{genfuj} 
  If $L$ is a $D$-big line bundle, then 
  \begin{equation}\label{eq:fujita} 
    \vol_{X|D}(L)=\langle L^{n-1}\rangle|_{\fD}. 
  \end{equation} 
\end{theo} 
This result is in fact established for a \emph{smooth}  
projective variety $X$, but it immediately extends to the case when  
$X$ is merely normal in view of the relation 
$\vol_{X|D}(L)=\vol_{X'|D'}(\mu^{\ast}L)$ if $\mu:X'\to X$  
is a blow-up and $D'$ is the strict transform of $D$  
(because $\mu$ has connected fibers, $X$ being normal). 
The theorem shows in particular that $\vol_{X|D}(L)$  
only depend on $c_1(L)\in N^1(X)$, and it follows that 
$\alpha\mapsto\langle\alpha^{n-1}\rangle|_{\fD}$  
is the (necessarily unique) extension of the restricted volume  
to a $(n-1)$-homogeneous and continuous function on the open cone  
of $D$-big classes. 
 
Using Theorem~\ref{charac}, we now extend this to  
arbitrary big line bundles. 
\begin{theo}\label{annul} 
  If $L$ is a big and $D$-psef line bundle 
  then~\eqref{eq:fujita} holds. 
  In particular, if $L$ is a big line bundle,  
  then $L$ is $D$-big iff   $\vol_{X|D}(L)>0$. 
\end{theo} 
As noted before, $L$ is $D$-big iff $D$ is not contained  
in the augmented base locus of $L$.  
We thus recover with a different proof a special case of the main  
result of~\cite{ELMNP3}.  
(The general case deals with irreducible components  
of the augmented base locus of any codimension.) 
\begin{remark}  
  If a line bundle $L$ is not $D$-psef,  
  then clearly $\vol_{X|D}(L)=0$, so the theorem shows  
  that $\vol_{X|D}(L)$ only depends on the numerical class  
  of the big line bundle $L$. It is worthwhile to note that  
  this fails for non-big line bundles in general. 
  For instance, if $D$ is (Cartier and) ample and $L$ is not big,  
  we have   $h^0(X,kL-D)=0$ for each $k$, and thus  
  $h^0(X|D,kL)=h^0(X,kL)$. Now if $X=C_1\times C_2$  
  is a product of two smooth curves and $L$ is the sum of the  
  pull-back of an ample bundle $L_1$ on $C_1$ and the  
  pull-back of a numerically trivial bundle $L_2$ on $C_2$,  
  then $\vol_{X|D}(L)=\limsup_{k\to\infty}\frac{1}{k}h^0(X,kL)=0$  
  when $L_2$ is not torsion in $\mathrm{Pic}^0(C_2)$, whereas 
  $\vol_{X|D}(L)=\limsup_{k\to\infty}\frac{1}{k}h^0(X,kL)=\deg L_1>0$  
  when $L_2$ is torsion. Note that $L$ is nef, so that we have  
  $\langle L\rangle|_{\fD}=L\cdot D=(D\cdot F_1)\deg L_1$,  
  where $F_1$ denotes the fiber of the projection $X\to C_1$.  
\end{remark}  
\begin{proof}[Proof of Theorem~\ref{annul}] 
  Let $L$ be an arbitrary $D$-psef line bundle. We claim that 
  \begin{equation}\label{e:mmm} 
    \langle L^{n-1}\rangle|_{\fD}\geq\vol_{X|D}(L). 
  \end{equation} 
  Indeed, if $A$ is an ample divisor on $X$, then clearly  
  $\vol_{X|D}(L+\ep A)\geq\vol_{X|D}(L)$ for every rational $\ep>0$.  
  Since $L+\ep A$ is $D$-big, the left-hand side coincides with  
  $\langle(L+\ep A)^{n-1}\rangle|_{\fD}$ by Theorem~\ref{genfuj}.  
  Thus~\eqref{e:mmm} follows by letting $\ep\to 0$,  
  by definition of $\langle L^{n-1}\rangle|_{\fD}$  
  for the $\fD$-psef class determined by $L$.  
   
  Now if $L$ is big and $D$-psef but not $D$-big,  
  then $\langle L^{n-1}\rangle|_{\fD}=0$ by Theorem~\ref{charac},  
  and thus~\eqref{e:mmm} gives 
  $\vol_{X|D}(L)=0$, which completes the proof. 
\end{proof} 
 
 
\subsection{Proof of Theorem~B} 
Let $L$ be a big line bundle on $X$. We also assume that the prime 
divisor $D$ is Cartier on $X$, which of course holds if $X$ is 
smooth. We have to show that $\vol_{X|D}(L)=\langle 
L^{n-1}\rangle\cdot D$. If $L$ is not $D$-big, then the left-hand side 
is zero by Theorem~\ref{annul}, and the right-hand side is zero by 
Lemma~\ref{orthobis}.  We can thus assume that $L$ is $D$-big. In that 
case, we have  
$\langle L^{n-1}\rangle\cdot D\geq\langle L^{n-1}\rangle|_{\fD} 
=\vol_{X|D}(L)$  
by Proposition~\ref{compar} and 
the generalized Fujita theorem. In order to prove the converse 
inequality, we rely on the following two simple remarks.  First, for 
every line bundle $M$ the kernel of the restriction map  
$H^0(X,M)\to H^0(D,M|_{D})$ is precisely $H^0(X,M-D)$, hence 
\begin{equation}\label{basic1} 
  h^0(X,M)-h^0(X,M-D)=h^0(X|D,M). 
\end{equation} 
Second, for any effective divisor $B$ whose support does not contain 
$D$, multiplication by the canonical section $\sigma$ of $H^0(X,B)$ 
yields injections  
$H^0(X,M-B) \hookrightarrow H^0(X,M)$ and 
$H^0(D,M-B) \hookrightarrow H^0(D,M)$.  
This implies the second basic relation:  
\begin{equation}\label{basic2} 
  h^0(X|D,M-B)\le h^0(X|D,M). 
\end{equation} 
Now fix an integer $k$, apply~\eqref{basic1} to  
$M=kL - jD$ for $j= 0, \dots, k$,  
and sum all these relations. This gives 
\begin{equation*} 
  h^0(X,kL) - h^0(X,k(L-D)) = \sum_{j=0}^{k-1} h^0(X|D, kL-jD). 
\end{equation*} 
Fix a sufficiently ample divisor $A$ not containing $D$ in its support 
such that $A-D$ is linearly equivalent to an effective divisor $B$ not 
containing $D$ in its support. By repeated use of~\eqref{basic2}, we 
get $h^0(X|D,kL-jD)\le h^0(X|D,kL-jD+jA)=h^0(X|D,kL+jB)\le 
h^0(X|D,kL+kB)$ for each $j\le k$. Hence 
\begin{equation*} 
  h^0(X,kL) - h^0(X,k(L-D)) \le k\,h^0(X|D, k(L+B)). 
\end{equation*} 
Dividing by $k^n$, and taking the limit when $k\to \infty$, we infer 
\begin{equation*} 
  \vol (L) - \vol( L-D) \le n\vol_{X|D}(L+B). 
\end{equation*} 
Replacing $L$ by $kL$, and expressing the (restricted) volumes as 
(restricted) positive intersection products, we finally conclude that 
\begin{equation*} 
  \frac1k\left(\langle L^n\rangle-\langle(L-\frac1k D)^n\rangle\right) 
  \le n\langle(L+\frac1k B)^{n-1}\rangle|_{\fD}. 
\end{equation*} 
When $k\to \infty$, the left-hand side tends to $n~\langle 
L^{n-1}\rangle \cdot D$ by Theorem~A, whereas the right-hand side 
converges to $n~\langle L^{n-1}\rangle|_{\fD}=n\vol_{X|D}(L)$ by 
continuity (Proposition~\ref{cont-rest}).  This proves the required 
inequality $\langle L^{n-1}\rangle\cdot D \leq \vol_{X|D}(L)$ and 
concludes the proof. 
 
\begin{example}[Surfaces] 
  In dimension $n=2$, a big class $\alpha\in \CNS(X)$ is not $D$-big 
  iff $P(\alpha)\cdot D=0$, where $P(\alpha)$ is the nef part of its 
  Zariski decomposition.  In particular, $\alpha$ is always $D$-big 
  when $D$ is ample (or even nef by the Hodge index theorem).  
  Now the continuous extension of $L\mapsto\vol_{X|D}(L)$  
  to the big cone of $\CNS(X)$, 
  to wit $\alpha\mapsto(P(\alpha)\cdot D)$, is not $\mathcal{C}^1$ in 
  general on the open cone of $D$-big classes, because $\alpha\mapsto 
  P(\alpha)$ is not $\mathcal{C}^1$ in general, as exemplified by the 
  blow-up of ${\bf P}^2$ at a point.  This means that the analogue of 
  Theorem~A fails for restricted volumes. 
\end{example} 
 

\end{document}